%% file: draft.tex
\documentclass[letterpaper, 10pt, conference]{ieeeconf}  

\IEEEoverridecommandlockouts                              

\makeatletter
\let\NAT@parse\undefined
\makeatother

\usepackage{times} 

\usepackage{cite}
\usepackage{amssymb}
\usepackage[tbtags]{amsmath}
\usepackage{algorithmic}
\usepackage{graphicx}
\usepackage{versions}

\let\labelindent\relax

\usepackage{amsthm}
\usepackage{thmtools}

\usepackage{calc}
\usepackage{ifthen}
\usepackage{afterpage}
\usepackage{luatex85, shellesc} 

\usepackage{stackengine} 
\usepackage{scalerel}   
\usepackage{float}

\usepackage[fleqn]{mathtools}
\usepackage[]{amsmath}
\PassOptionsToPackage{fleqn}{amsmath}
\usepackage{interval}
\usepackage{physics}

\usepackage{bookmark}

\usepackage[hang]{footmisc}

\usepackage[inline]{enumitem}
\usepackage{xspace}
\usepackage{siunitx}

\usepackage{tabularx}
\usepackage{ltablex}
\usepackage{multicol}
\keepXColumns

\usepackage{makecell}

\usepackage{booktabs}

\usepackage{tikz}               
\usepackage{tikz-3dplot}
\usepackage{pgfplots}
\usepackage{tikzscale}
\usepackage{pgf}
\usepackage{pgfplotstable}

\usepgfplotslibrary{groupplots}
\usepgfplotslibrary{external}

\usetikzlibrary{calc}
\usetikzlibrary{arrows.meta}
\usetikzlibrary{backgrounds}
\usetikzlibrary{trees}
\usetikzlibrary{positioning}
\usetikzlibrary{shapes.multipart}

\usetikzlibrary{external}
\usepackage{hyperref}
\usepackage{cleveref}

\input{helper_setup.tex}

\title{\LARGE \bf
A Note on Nussbaum-type Control and Lie-bracket Approximation
}

\author{Marc Weber$^{1}$\!\!,\,\, Christian Ebenbauer$^{1}$ and Bahman Gharesifard$^{2}$
\thanks{$^{1}$\,Chair of Intelligent Control Systems, RWTH Aachen University, Germany
{\tt\small \href{mailto:marc.weber@ic.rwth-aachen.de}{marc.weber@ic.rwth-aachen.de}}, {\tt\small \href{mailto:christian.ebenbauer@ic.rwth-aachen.de}{christian.ebenbauer@ic.rwth-aachen.de}}}%
\thanks{$^{2}$\,Department of Electrical and Computer Engineering, University of California, Los Angeles
        {\tt\small \href{mailto:gharesifard@ucla.edu}{gharesifard@ucla.edu}}}%
}

\begin{document}

\maketitle
\thispagestyle{empty}
\pagestyle{empty}

\begin{abstract}
    In this paper, we propose an adaptive control law for completely unknown scalar linear systems based on  Lie-bracket approximation methods. We investigate stability and convergence properties for the resulting Lie-bracket system, compare our proposal with existing Nussbaum-type solutions and demonstrate our results with an example.
    Even though we prove global stability properties of the Lie-bracket system, the stability properties of the proposed dynamics remain open, making the proposed control law an object of further studies.
    We elaborate the difficulties of establishing stability results by investigating connections to partial stability as well as studying the corresponding Chen-Fliess expansion.
\end{abstract}

\section{PROBLEM FORMULATION}
\label{sec:problem-fomulation}
Consider the linear time-invariant plant
\begin{align}
    \label{eqn:lti-problem}
    \dot y(t) & = ay(t) + bu(t)
\end{align}
with unknown $a\in \R$ and unknown $b \in \Rnot{0}$. Our goal is to find a preferrably smooth control law, which does not depend on the knowledge of $a$ and $b$, such that for all $y_0 \in \R$ and for all $t_0 \in \R$, the solution $y(t; t_0, y_0)$ converges to zero, i.e. $\lim_{t \to \infty} y(t; t_0, x_0) = 0$, and all control signals are at least bounded.

In 1983, Morse conjectured, that no universal stabilizer as sought in the problem formulation can exist for linear systems, which is of the same dynamic order as the plant \eqref{eqn:lti-problem} and rational or polynomial in controller and plant states \cite{nussbaum1983some}. Nussbaum proved in 1984 that Morse's conjecture was indeed true and gave rise to a whole class of possible control algorithms \cite{nussbaum1983some}.
A multivariable version was given the same year by Byrnes and Willems \cite{byrnes1984adaptive}.
In the following years, the topic was researched in \cite{heymann1985remarks}, \cite{maartensson1985order}, \cite{ilchmann1987high}, \cite{owens1987positive} and \cite{morse1984adaptive}. The report \cite{helmke1988stability} summarises several contributions.

After decades of little attention, Nussbaum's design has been incorporated in observers \cite{talkhoncheh2017observer} and model reference adaptive controllers \cite{gerasimov2018relaxing} recently. While different definitions of Nussbaum-type functions and control designs exist \cite{chen2019nussbaum}\cite{ilchmann1993adaptive}, we have selected the following one from \cite{ilchmann1993adaptive}.

A piecewise right continuous map $h : \R \to \R$ is said to be of \emph{Nussbaum-type}, if for some $k_0 \in \R$
\begin{subequations}
    \label{eqn:nussbaum-condition-ilchmann}
    \begin{align}
        \sup_{k > k_0} \frac 1 {k-k_0} \int_{k_0}^k h(s) s \dx[s] & = + \infty
        \intertext{and simultaneously}
        \inf_{k > k_0} \frac 1 {k-k_0} \int_{k_0}^k h(s) s \dx[s] & = - \infty.
    \end{align}
\end{subequations}
Assume $h$ is of Nussbaum-type. Then for all $a \in \R$, $b \in \Rnot{0}$, $y_0 \in \R$, $k_0 \in \R$ and $t_0 \in \R$, the solutions $y(t; t_0, y_0)$ and $k(t; t_0, k_0)$ of \eqref{eqn:lti-problem} with the control law
\begin{subequations}
    \label{eqn:nussbaum}
    \begin{align}
        \dot{k}(t)    & = y^2(t)    \label{eqn:nussbaum:adaptive-law}            \\
        u(x(t), k(t)) & = h(k(t))\, k(t)\, y(t) \label{eqn:nussbaum:control-law}
    \end{align}
\end{subequations}
are defined for all $t \geq t_0$, $\lim_{t \to \infty} y(t; t_0, y_0) = 0$ and $\lim_{t\to \infty} k(t; t_0, k_0)$ exists and is finite \cite{morse1984adaptive}.

One commonly used Nussbaum-type function is {$h(s) = s \cos(s)$}.

A contemporary approach called funnel control, as summarized in \cite{berger2021funnel}, generalizes the concept of Nussbaum functions and is also able to handle unknown control directions for our problem \eqref{eqn:lti-problem}. Funnel control uses a prescribed funnel boundary to have a designer-specified convergence rate, and in turn, this funnel boundary introduces an explicit time-dependency into the control law, which renders the closed loop non-autonomous.

For the more common setting in adaptive control, where the sign of $b$ is known, the Willems-Byrnes control algorithm \cite{willems1984global} provides a very intuitive solution. Consider
\begin{subequations}
    \label{eqn:willems-byrnes-control}
    \begin{align}
        \dot{k}(t)    & = \operatorname{sign}(b)\, y^2(t) \quad \text{and} \\
        u(y(t), k(t)) & = -k(t)\,y(t).
    \end{align}
\end{subequations}
For all unknown $a \in \R$, $b \in \Rnot{0}$ with known $\operatorname{sign}(b)$, $y_0 \in \R$, $k_0 \in \R$ and $t_0 \in \R$, the solutions $y(t; t_0, y_0)$ and $k(t; t_0, k_0)$ are defined for all $t \geq t_0$, $\lim_{t \to \infty} y(t; t_0, y_0) = 0$ and $\lim_{t \to \infty} k(t; t_0, k_0)$ exists and is finite.

While a lot of effort has previously been spent on researching and advancing the class of Nussbaum-type control and adaptive laws, we hit the topic from the perspective of online optimization by using Lie-bracket approximation of trajectories \cite{durr2013lie} \cite{morse1984adaptive} \cite{labar2019thesis}.
In \cite{6760498}, for a system of the form $\dot{y}= f(y) +g(y)u$, a control law was proposed, which approximated the damping control $u= -\alpha \LieDeriv{y}{g}{V}$ through Lie-bracket approximation without the need to know the control vector field $g$, however a lower bound on the gain $\alpha$ was still required.
In a conclusive comment in \cite{scheinker2012minimum}, the authors briefly touch on a similar universal controller design but note that the particular exemplified design, which is an adaptive extension to the main results in \cite{scheinker2012minimum}, only can be shown to be stable in a semi-global practical sense and thus it would be expectected for an adapted gain $k(t)$ to drift towards infinity.

The idea of this work is to extend the intuitive Willems-Byrnes control and update law such that the sign of $b$ is not needed for implementation and without a need to introducing Nussbaum-type functions.
Towards that end, it is necessary to expose the input gain~$b$ through some other means. As we demonstrate, this can be achieved using Lie-bracket approximation.

This following paragraph on Lie-bracket approximation sums up the required results from \cite{labar2019thesis}. Definitions and further results discussed are taken from \cite{durr2013lie} and  \cite{grushkovskaya2018class}.
Consider the input-affine system
\begin{align}
    \label{eqn:labar:system}
    \dot{x}(t) & = f_0(t, x(t)) + \sum_{i=1}^l \omega^{p_i} f_i(t, x(t))\, u_i(k_i \omega t),
\end{align}
with $x(t) \in \R^n$ the state vector of the system, $u_i(t) \in \R$ the control inputs, $\omega \in \R_{>0}$, $k_i \in \Q$ and $p_i \in \interval[open]{0}{1}$.

Under some technical assumptions (\Cref{assumption:ui,assumption:fi,assumption:pi}), a so-called \emph{Lie-bracket system} may be associated with \eqref{eqn:labar:system}, namely
\begin{align}
    \label{eqn:labar:lbs}
    \dot{\bar{x}}(t) & = f_0(t, \bar{x}(t)) + \lim_{\omega \to \infty}\sum_{\substack{1 \leq i < l \\ i < j \leq l}} [f_i, f_j](t, \bar{x}(t)) \, \gamma_{ij}(\omega),
\end{align}
where we introduced
\begin{align}
    \label{eqn:labar:gamma}
    \gamma_{ij}(\omega) & \coloneqq \frac{\omega^{p_i + p_j}}{T} \int_0^T \int_0^\theta u_j(k_j \omega \theta) \, u_i(k_i \omega \tau) \dx[\tau]\dx[\theta],
\end{align}
with\footnote{The function $\operatorname{LCM}(\cdot)$ denotes the least common multiple operator.} $T = 2 \pi \omega^{-1} \operatorname{LCM}(k_1^{-1}, k_2^{-1}, \ldots, k_l^{-1})$.

Labar's results in \cite{labar2019thesis} allow us then to approximate solutions to \eqref{eqn:labar:system}  by solutions to \eqref{eqn:labar:lbs} with arbitrary accuracy over any finite timespan by choosing any suitably high $\omega$. The exact wording and technical premises are given in the Appendix as \Cref{prop:labar-1.1}.

\section{RESULTS}
\label{sec:results}

\subsection{Algorithm}
\label{sec:algorithm}
For the system \eqref{eqn:lti-problem} define the update and control law as
\begin{subequations}
    \label{eqn:mycontrol}
    \begin{align}
        \dot{k}(t)    & = y^2(t) \sqrt{\omega} \,\cos(\omega t) \quad \text{and} \label{eqn:mycontrol:update-law} \\
        u(x(t), k(t)) & = - k(t)\,y(t) - y(t) \sqrt{\omega}\, \sin(\omega t)
        \label{eqn:mycontrol:control-law}
    \end{align}
\end{subequations}
respectively with $\omega \in \interval[open]{\omega^*}{\infty}$ and $\omega^* \in \R_{> 0}$ sufficiently large. The set $X_e \coloneqq \{0\} \times \R$ is the equilibrium set for the closed-loop system of \eqref{eqn:lti-problem}, \eqref{eqn:mycontrol}. Notice that $X_e$ is not bounded.
We will proceed to motivate this design step by step.

The reader can observe that the control law \eqref{eqn:mycontrol} along with \eqref{eqn:lti-problem} can be realized as an instance of the control-affine system \eqref{eqn:labar:lbs}.
Consequently we can apply Lie-bracket approximation and calculate the associated Lie-bracket system.

\begin{proposition}[Lie-bracket system]
    \label{prop:lbs:cartesian}
    The Lie-bracket system associated with \eqref{eqn:lti-problem} and \eqref{eqn:mycontrol} is given as
    \begin{subequations}
        \label{eqn:dynamics-lbs-cartesian}
        \begin{align}
            \dot{\bar{y}}(t) & = (a - b \bar{k}(t))\,\bar{y}(t)                       \\
            \dot{\bar{k}}(t) & = b \bar{y}^2(t). \label{eqn:dynamics-lbs-cartesian-k}
        \end{align}
    \end{subequations}
\end{proposition}
The \hyperref[proof:lbs:cartesian]{proof} is presented in the Appendix. The Lie-bracket system \eqref{eqn:dynamics-lbs-cartesian} exposes $b$ into the dynamics of the update law, quite similar to \eqref{eqn:willems-byrnes-control}, while the control and update laws \eqref{eqn:mycontrol} can be implemented without the knowledge of $b$.

Lie-bracket approximation results assert that for any fixed window of time and for any arbitrary prescribed approximation error, we can find $\omega^*$ such that the trajectory of \eqref{eqn:lti-problem} and \eqref{eqn:mycontrol} follows closely the trajectory of \eqref{eqn:dynamics-lbs-cartesian} with $y_0 = \bar{y}_0$ and $k_0 = \bar{k}_0$ for all $\omega > \omega^*$.
Our next result aims at stability properties of this Lie-bracket system.

\begin{theorem}[Lie-bracket system analysis]
    \label{thm:lbs-system-analysis-cartesian}

    Consider the Lie-bracket system \eqref{eqn:dynamics-lbs-cartesian}. For any $p \in \R_{\geq0}$, let ${c_p = \frac{a + p}{b}}$ and $d_p(\bar{y}_0, \bar{k}_0) = \frac 1 2 ( \bar{y}_0^2 + (\bar{k}_0 - c_p)^2)$. Denote ${\bar{x}(t) = \begin{bmatrix}\bar{y}(t) & \bar{k}(t)\end{bmatrix}^T}$. Then,
    \begin{enumerate}[label=\roman*., ref=(\roman*)]
        \item \label[statement]{itm:lbs-1}the set $X_e \coloneqq \{0\}\times \R$ is the equilibrium set for \eqref{eqn:dynamics-lbs-cartesian};
        \item \label[statement]{itm:lbs-2} any point $\bar{x}_p = (0, c_p) \in X_e$ is globally stable, but not asymptotically stable;
        \item \label[statement]{itm:lbs-3}the solutions $\bar{y}(t; t_0, y_0)$ and $\bar{k}(t; t_0, k_0)$ exist for all $t \geq t_0$ and are unique;
        \item \label[statement]{itm:lbs-4}$\bar{x}(t; t_0, x_0) \in \overline{\mathcal{B}_{d_p}}((0, c_p))$ for all $x_0 \in \R^2$; and
        \item \label[statement]{itm:lbs-5} $\lim_{t \to \infty} \bar{x}(t; t_0, x_0) = (0, c_0 + \operatorname{sign}(b) \, d_0)$ for all $x_0 \in \R^2 \setminus X_e$.
    \end{enumerate}
\end{theorem}

The \hyperref[proof:lbs:properties]{proof} is given in the Appendix.
For a similar system to \eqref{eqn:dynamics-lbs-cartesian}, boundedness of all closed-loop signals and convergence of $\bar{x}$ has been proven using the same inference chain in the context of adaptive output-feedback control, e.g. in \cite{fradkov2013nonlinear} under the restriction $\omega(t) \equiv 0$. However, the mentioned proof requires knowledge of the high frequency gain implicitely.

Our proof notably illustrates that the orbits of the Lie-bracket system form circular arcs. Therefore, it is only natural to try and express the dynamics in polar coordinates.

\begin{proposition}[Transformation to polar coordinates]%
    \label{prop:polar-coordinate-transformation}%
    Let $c_0 \in \R$ be arbitrary but fixed.
    Consider the transformation
    \begin{align}
        r(y, k)       & = \sqrt{y^2 + (k- c_0)^2},  \label{eqn:transform-r}           \\
        \varphi(y, k) & = \begin{cases}
                              \arcsin \frac{k - c_0}{r(y, k)}       & \text{if } y \geq 0 \\
                              \pi - \arcsin \frac{k - c_0}{r(y, k)} & \text{if } y < 0
                          \end{cases}
        \label{eqn:transform-phi}
    \end{align}
    with inverse transformation
    \begin{align}
        y(r, \varphi) & = r \cos \varphi        \\
        k(r, \varphi) & = r \sin \varphi + c_0.
    \end{align}
    In the $(r, \varphi)$-coordinates, \eqref{eqn:lti-problem} with \eqref{eqn:mycontrol} can be expressed as
    \begin{subequations}
        \begin{align}
            \phantom{\dot{r}(t) =}
             &
            \begin{aligned}[b]
                \mathllap{\dot{r}(t) =}- \,  br^2(t) \sin \varphi(t) \cos^2 \varphi(t) &                                     \\
                -\,  b r(t) \cos^2 \varphi(t)                                          & \cdot\sqrt{\omega} \, u_1(\omega t) \\
                + \, r^2(t) \sin \varphi(t) \cos^2 \varphi(t)                          & \cdot\sqrt{\omega} \, u_2(\omega t)
            \end{aligned}
            \\
            \phantom{\dot{\varphi}(t)  =}
             &
            \begin{aligned}[b]%
                \mathllap{\dot{\varphi}(t)  = + } \, b r(t) \sin^2 \varphi(t) \cos \varphi(t) &                                      \\
                +\, b \sin \varphi(t) \cos \varphi(t)                                         & \cdot\sqrt{\omega} \, u_1(\omega t)  \\
                +\, r(t) \cos^3 \varphi(t)                                                    & \cdot\sqrt{\omega} \, u_2(\omega t).
            \end{aligned}
        \end{align}
        \label{eqn:mycontrol-closedloop-polar}
    \end{subequations}
    Likewise, the Lie-bracket system from equation~\eqref{eqn:dynamics-lbs-cartesian} can be expressed in $(r, \varphi)$-coordinates as
    \begin{subequations}
        \label{eqn:dynamics-lbs-polar}
        \begin{align}
            \dot{\bar{r}}(t)       & = 0                                   \\
            \dot{\bar{\varphi}}(t) & = b \bar{r}(t) \cos \bar{\varphi}(t).
        \end{align}
    \end{subequations}%
\end{proposition}%
\begin{conferenceversion}%
    The proof of \Cref{prop:polar-coordinate-transformation} is given in the extended version~\cite{weberebenbauergharesifard2021}.%
\end{conferenceversion}%
\begin{extendedversion}%
    The \hyperref[proof:polar-coordinates]{proof} is given in the Appendix.%
\end{extendedversion}%
In polar coordinates, the Lie-bracket system shows more articulately, how the orbits of \eqref{eqn:dynamics-lbs-cartesian} respectively \eqref{eqn:dynamics-lbs-polar} describe circular arcs with constant polar radius $r$, which has already been established in the proof of \Cref{thm:lbs-system-analysis-cartesian}.

\subsection{Discussion}
\label{sec:conclusions-remarks-discussion}

We now discuss in details the difficulties that arise in establishing (practical) stability properties for the original dynamics \eqref{eqn:lti-problem} and \eqref{eqn:mycontrol}.

While we have shown global stability and convergence properties of the Lie-bracket system, one cannot infer any (practical) stability result on the system \eqref{eqn:lti-problem},\eqref{eqn:mycontrol} using existing Lie-bracket approximation results. The reason for this is that current systematic results on Lie-bracket approximation rely on uniform asymptotic stability properties of the Lie-bracket system. Moreover, the stability result of the Lie-bracket system cannot be strengthened to asymptotic stability as we have shown in \Cref{thm:lbs-system-analysis-cartesian}. This being said, \Cref{prop:labar-1.1} allows us to make statements about the behaviour of the controlled system \eqref{eqn:lti-problem}, \eqref{eqn:mycontrol} over finite times when the initial conditions are selected from a bounded region. However, due to the lack of uniform asymptotic stability of any equilibrium point of the Lie-bracket system, established results such as the ones in \cite{labar2019thesis} cannot be applied directly.

We have also investigated the partial stability notion in \cite{vorotnikov2005partial} to show uniform $y$-stability of the whole equilibrium set $X_e$.

One can deduce from \Cref{thm:lbs-system-analysis-cartesian} that indeed the point $(0, c_0)$ is globally $y$-stable for the Lie-bracket system, a natural consequence from its global stability in the sense of Lyapunov, and moreover, the points $(0, c_p)$ for $p > 0$ are globally asymptotically y-stable according to the direct application of \cite[Definition 1.(2,6)]{vorotnikov2005partial}. The asymptotic $y$-stability is naturally uniform in $t_0$ for our Lie-bracket system. However, the convergence of $y(t; t_0, x_0)$ is not uniform in the initial conditions $x_0$. This has been shown in the recent work \cite{orlowski2020}, where a special case of our Lie-bracket system (taking $a=0$, $b=1$ and restricting the domain to $D=\{x: \|x-x^*\|\leq 1\}$) has been shown to lack the uniformity part of uniform asymptotic y-stability.
As also detailed in \cite{orlowski2020}, this means that the overshoot, which in our design is inevitable as part of the adaptation of $k$, can happen arbitrarily late.

Our extensive numerical experiments suggest that the trajectories of \eqref{eqn:lti-problem},\eqref{eqn:mycontrol} are bounded and converge under the condition that a sufficiently high $\omega$ has been chosen. However, we observed that for some large initial conditions the simulations are subject to numerical instabilities, and hence it is unclear, if this is attributed to an intrinsic instability of the dynamics.

We wish to emphasize that our simple design gives rise to a large class of control laws with the same, or at least similar, Lie-bracket systems. They all share the exploitation of Lie-bracket approximation to bypass the unknown parameters $a$ and $b$. One other design of this class are the control and adaptive laws
\begin{subequations}
    \label{eqn:mycontrol-swapped}
    \begin{align}
        \dot{k}       & = y(t) \sqrt{\omega} \,\cos(\omega t) \quad \text{and} \label{eqn:mycontrol-swapped:update-law} \\
        u(y(t), k(t)) & = - k(t)\,y(t) - 2y^2(t) \sqrt{\omega}\, \sin(\omega t).
        \label{eqn:mycontrol-swapped:control-law}
    \end{align}
\end{subequations}
Notice the components of the vectorfields in \eqref{eqn:mycontrol-swapped} have been swapped and an additional factor of \num{2} has been introduced in order to provide the identical associated Lie-bracket system \eqref{prop:lbs:cartesian}.
We observe in \Cref{subsec:swap} that even though both designs share the same Lie-bracket system, their qualitative and quantitative behaviours for the same $\omega$ differ noticeably. In part, this difference can be attributed to the additional factor introduced in \eqref{eqn:mycontrol-swapped:control-law}, which in turn will lead to a change in the bound $M$ of \Cref{assumption:fi-2}.

On the other hand, comparison with Nussbaum-type control laws in \cref{subsec:example-comparison-with-nussbaum} will show, that our design seems to be numerically more reliable. While convergence time is in general (that is, one can always find initial conditions, such that one design or the other shines) faster for Nussbaum-type control laws, our design exhibits a better transient behaviour with less overshoot of the controlled state $y$. Both control and adaptive laws can be further tuned using separate gains.

We leave comparison with a funnel control law for future work, as we are indecisive on a suitable choice of boundary for a meaningful comparison.

We also investigated the functional expansion using a Chen-Fliess series of the solution $x(t; t_0, x_0)$ similar to the functional expansion in \cite{labar2019thesis}, see \Cref{thm:cf-expansion}, to study higher-order Lie-bracket approximations. While some simplifications can be made, it appears to be difficult to bound the occurring terms in a significantly tighter way. We could not make one resulting sufficient choice of $\omega^*$ global, in fact the sufficient choice of $\omega^*$ depends (for each iterated application) on the respective initial condition.
\begin{conferenceversion}
    This being said, our calculations, which are part of the extended version on arXiv \cite{weberebenbauergharesifard2021}, show that only few iterated Lie derivatives contribute to the output, and a more careful study of the pattern present in the terms of this series may lead to conclusive results. However, the power of the contributing states in some cases is higher than the power of the time period, which seems to indicate that for any fixed time period or $\omega$, some initial conditions will exist for which one of the states might grow unbounded.
\end{conferenceversion}
\begin{extendedversion}
    Our calculations, which follow in \cref{subsec:chen-fliess-expansion}, show that only few iterated Lie derivatives contribute to the output. However, the power of the contributing states in some cases is higher than the power of the time period, which seems to indicate that for any fixed time period or $\omega$, some initial conditions will exist for which one of the states might grow unbounded.
\end{extendedversion}
Due to these obstacles, we believe that the second order system \eqref{eqn:lti-problem},\eqref{eqn:mycontrol} defines a challenging benchmark example for future research in stability and Lie-bracket approximation theory.

\section{NUMERICAL EXPERIMENTS}

We demonstrate the application of our design using numerical experiments.

\subsection{Demonstration}
\label{subsec:example-demonstration}
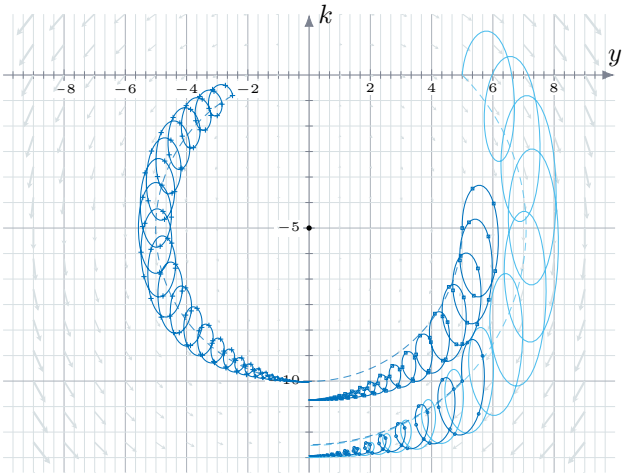
\begin{figure}[tb]%
    \tikzsetnextfilename{figure_orbits}%
    \centering%
    \input{figures/figure_orbits.tikz}%
    \caption{Simulation results: Graphs show orbits of \eqref{eqn:lti-problem},\eqref{eqn:mycontrol} in solid and orbits of \eqref{eqn:dynamics-lbs-cartesian} in dashed over a quiver plot illustrating the vectorfield of the Lie-bracket system \eqref{eqn:dynamics-lbs-cartesian}. }%
    \label{fig:example}
\end{figure}

Let $a=10$ and $b=-2$ as well as $\omega=400$. We simulate the system \eqref{eqn:lti-problem},\eqref{eqn:mycontrol} using an {\itshape ode1}-solver with stepwidth $h= \frac 1 {40} T = \frac{2 \pi}{40 }\omega^{-1}$. The parameters are arbitrary choices and a compromise between enough resolution to see details and reasonable computational effort. The results are depicted in \Cref{fig:example}.

\subsection{Comparison}
\label{subsec:example-comparison-with-nussbaum}

To illustrate possible benefits, we also simulate \eqref{eqn:lti-problem} with a Nussbaum-type control and adaptation law \eqref{eqn:nussbaum} and compare it with our proposed alternative.
\begin{figure}[tb]%
    \tikzsetnextfilename{figure_trajectories}%
    \centering%
    \includegraphics[width=\linewidth]{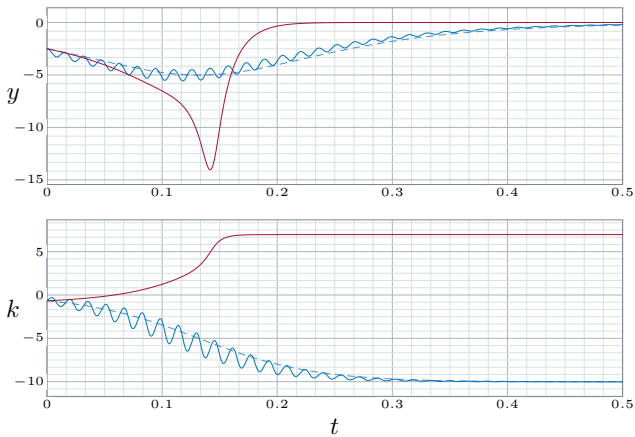}
    \caption{Graphs show trajectories of our proposed control \eqref{eqn:lti-problem},\eqref{eqn:mycontrol} as \ref{plot:ex2:es:y:mw} and its associated Lie-bracket system \eqref{eqn:dynamics-lbs-cartesian} as \ref{plot:ex2:lbs:y:mw} in comparison to trajectories of a Nussbaum-type control \eqref{eqn:nussbaum} with $h(s) = s\cos(s)$ \eqref{eqn:lti-problem} as \ref{plot:ex2:y:nblm}.}%
    \label{fig:comparison-nussbaum}%
\end{figure}
We use the Nussbaum-type function $h(s) = s\cos s$ and an {\itshape ode1}-solver, however, we had to decrease the stepwidth to $h=\num{1e-4}$ in order for the problem to remain numerically stable. The comparison is displayed in \Cref{fig:comparison-nussbaum}.

\subsection{Different algorithm}
\label{subsec:swap}
\begin{figure}[hbt]%
    \tikzsetnextfilename{figure_orbits_swapped}%
    \centering%
    \input{figures/figure_orbits_swapped.tikz}%
    \caption{Simulation results for \Cref{subsec:swap}: Graphs show orbit of \eqref{eqn:lti-problem},\eqref{eqn:mycontrol} as \ref{plot:ex-swapped:es}, orbit of \eqref{eqn:lti-problem},\eqref{eqn:mycontrol-swapped} as \ref{plot:ex-swapped:es-swapped} and orbit of their common Lie-bracket system \eqref{eqn:dynamics-lbs-cartesian} as \ref{plot:ex-swapped:lbs} over a quiver plot illustrating the vector field of the common Lie-bracket system \eqref{eqn:dynamics-lbs-cartesian}. }%
    \label{fig:example-swapped}%
\end{figure}
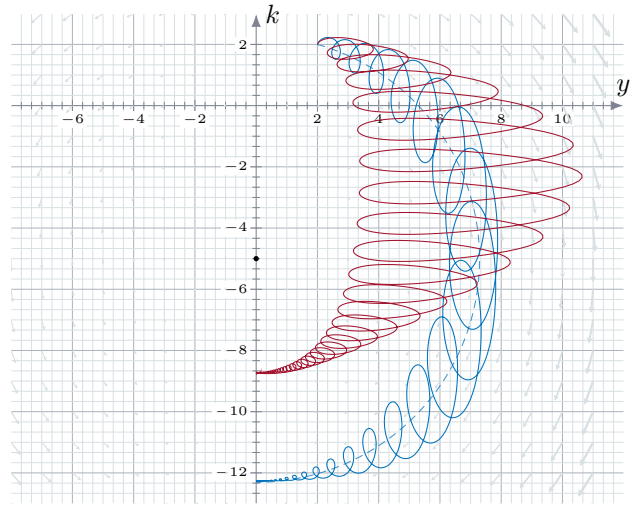

We simulate the system \eqref{eqn:lti-problem} and \eqref{eqn:mycontrol-swapped} and compare it for one initial condition with our original design \eqref{eqn:mycontrol}.  We reused the parameters from \Cref{subsec:example-demonstration}. The comparison is shown in \Cref{fig:example-swapped}.

\subsection{Chen-Fliess expansion}
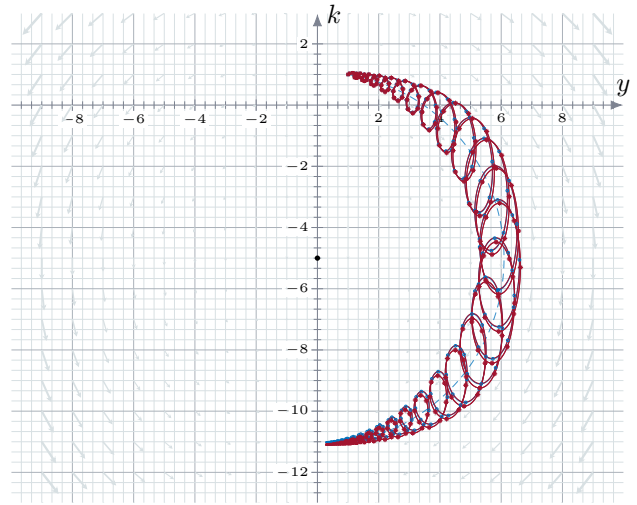
\begin{figure}[hbt]%
    \tikzsetnextfilename{figure_orbits_cf}%
    \centering%
    \input{figures/figure_orbits_cf.tikz}%
    \caption{ Lie bracket system as \ref{plot:ex3:lbs}, our proposed control as \ref{plot:ex3:es} and Chen-Fliess expansion integration schem of our proposed control of order $d=0$ as \ref{plot:ex3:cf0}, order $d=1$ as \ref{plot:ex3:cf1} and order $d=2$ as \ref{plot:ex3:cf2} for a common initial condition.}%
    \label{fig:example-cf}%
\end{figure}

We use the Chen-Fliess expansion of orders $d=0, 1, 2$ as a numerical integration scheme for a common initial condition and compare the resulting orbits with one orbit from \Cref{subsec:example-demonstration}. For numerically solving the ode, we used the same settings as in \Cref{subsec:example-demonstration}, and for the Chen-Fliess expansion integration scheme we used a stepwidth $T = \frac{2\pi}{64}\omega^{-1} \approx \num{2.5e-4}$. The comparison is shown in \Cref{fig:example-cf}.

\section*{ACKNOWLEDGMENT}

Bahman Gharesifard is thankful to the Alexander von Humboldt foundation for their support.


\bibliography{sources}
\bibliographystyle{IEEEtran}


\appendices

\section{Tools}
\label{sec:appendix:tools}

\subsection{Lie-bracket Averaging}
We recall some results  for Lie-bracket averaging following \cite{labar2019thesis}. Consider the input-affine control system \eqref{eqn:labar:system}.
\begin{assumption}[Conditions on  \texorpdfstring{$u_i$}{ui}]
    \label{assumption:ui}
    Assume that for all $i \in \{1, \ldots, l\}$ the following holds:
    \begin{assenum}
        \item \label{assumption:ui-1}
        $u_i: \R \to \R$ is a measurable function such that $\sup_{t\in \R} |u_i(t)| \leq 1$,
        \item \label{assumption:ui-2}
        $u_i(t)$ is $2\pi$-periodic, i.e. $u_i(t) = u_i(t+2\pi)$ for all $t \in \R$,
        \item  \label{assumption:ui-3}
        $u_i(t)$ has zero mean over one period, i.e.\! ${\int_0^{2\pi}\!\!\!\! u_i(t) \dx[t] = 0}$.
    \end{assenum}
\end{assumption}

\begin{assumption}[Conditions on \texorpdfstring{$f_i$}{fi}]
    \label{assumption:fi}
    We assume that there exists an open convex set $\mathcal{R} \subseteq \R^n$, such that for all  $i \in \{0, 1, \ldots, l\}$ and all $j \in \{1, \ldots, l\}$ the followings hold:

    \begin{assenum}
        \item \label{assumption:fi-1} $f_i : \R \times \mathcal{R} \to \R^n$ and $\LieDeriv{f_i}{f_j} : \R \times \mathcal{R} \to \R^n$ are differentiable vector fields,
        \item \label{assumption:fi-2} for every compact set $\mathcal{A} \subset \mathcal{R}$, there exists an $M \in \R_{> 0}$, such that for all $x \in \mathcal{A}$ and $t \in \R$, the following bounds apply:
        \begin{multicols}{2}
            \begin{enumerate}[label=\roman*., ref=(\roman*)]
                \item $\norm{f_i(t, x)} \leq M$,
                \item $\norm{\mathcal{D}_t f_i(t, x)} \leq M$,
                \item $\norm{\mathcal{D}_x f_i(t, x)} \leq M$,
                \item $\norm{\mathcal{D}_t \LieDeriv[t, x]{f_i}{f_j}} \leq M$, and
                \item  $\norm{\mathcal{D}_x \LieDeriv[t, x]{f_i}{f_j}} \leq M$.
            \end{enumerate}
        \end{multicols}
    \end{assenum}

\end{assumption}

\begin{assumption}[First order conditions]
    \label{assumption:pi}
    Let $\mathcal{R} \subseteq \R^n$ come from \Cref{assumption:fi}. We assume that for all $i \in \{1, \ldots, l\}$ the vector fields and dithers satisfy the following conditions:
    \begin{assenum}
        \item \label{assumption:pi-1} If $p_i + p_j > 1$, then
        \begin{enumerate}[label=\roman*., ref=(\roman*)]
            \item $\int_0^{T_{ij}} \int_0^{\theta} u_j(k_j \omega \theta) u_i(k_i \omega \tau)\dx[\tau] \dx[\theta] = 0$, with $T_{ij}$ a common period to $u_i(k_i \omega t)$ and $u_j(k_j \omega t)$, or
            \item $\LieB[t, x]{f_i}{f_j} = 0$ for all $x\in \mathcal{R}$ and all $t \in \R$,
        \end{enumerate}

        \item \label{assumption:pi-2} If $p_i + p_j + p_m \geq 2$, then $\LieDeriv[t,x]{\LieDeriv{f_m}{f_j}}{f_i} = 0$ for all $x \in \mathcal{R}$ and all $t \in \R$.
    \end{assenum}
\end{assumption}

\begin{proposition}[Lie-bracket averaging on finite time intervals]
    \label{prop:labar-1.1}
    Suppose that \Cref{assumption:ui,assumption:fi,assumption:pi} hold and let $\mathcal{R} \subseteq \R^n$ be as in \Cref{assumption:fi}.
    Let $t_f \in \R_{\geq 0}$ and $\mathcal{V}, \mathcal{W} \subset \mathcal{R}$ be compact and convex compact respectively, such that for every $t_0 \in \R$ and $\bar{x}_0 \in \mathcal{V}$, the trajectory of system~\eqref{eqn:labar:lbs} through $\bar x (t_0) = \bar{x}_0$ exists, is unique and satisfies $\bar x(t) \in \mathcal{W}$ for all $t \in \interval{t_0}{t_0 + t_f}$.
    Then, for every $D \in \R_{> 0}$, there exists an $\omega^* \in \R_{>0}$ such that for all $\omega \in \interval[open]{\omega^*}{\infty}$, for all $t_0 \in \R$ and for all $x_0 \in \mathcal{V}$, the trajectories of systems~\eqref{eqn:labar:system} and \eqref{eqn:labar:lbs} through $x(t_0) = x_0$  and $\bar{x}(t_0) = x_0$ satisfy
    \begin{align}
        \norm{x(t; t_0, x_0) - \bar{x}(t; t_0, x_0)} < D\quad \forall t \in \interval{t_0}{t_0 + t_f}.
    \end{align}
\end{proposition}

\subsection{Chen-Fliess expansion}
\label{subsec:chen-fliess-expansion}
The following theorem has been adapted from \cite[Theorem 1]{qiao2010fliess}. Consider an affine nonlinear system
\vspace{-6pt}
\begin{subequations}
    \label{eqn:cf-system-abstract}
    \begin{align}
        \dot{x}(t) & = f_0(x) + \sum_{i=1}^m f_i(x(t)) \, u_i(t) \\
        y(t)       & = h(x(t)) = \begin{bmatrix}
                                     h_1(x(t)) & \dots & h_p(x(t))
                                 \end{bmatrix}^T,
    \end{align}
\end{subequations}
where $x(t) \in \R^n$, $u_i(t) \in \R$, $y(t) \in \R^p$, $f_0$ and $f_i, i = 1 \dots, m$ are smooth vector fields and $h_j, j=1, \dots, p$ are smooth functions.

\begin{theorem}
    \label{thm:cf-expansion}
    The input-output response of system \eqref{eqn:cf-system-abstract} can be expressed for $j=1\dots p$ as
    \begin{align}
        \label{eqn:cf-input-output-response}
        y_j(t)  =
        \begin{split}
            &\,h_j(x_0)                                                                                                 \\
            &+ \sum_{d=0}^{\infty} \sum_{i_d, ..., i_0 = 0}^{m} \LieDeriv{f_{i_0}}{\dots\LieDeriv{f_{i_k}}{h_j(x_0)}}
            \int_{t_0}^t \dx[\xi]_{i_d} \dots \dx[\xi]_{i_0}
        \end{split}
    \end{align}
    where the symbol $\int_{t_0}^t \dx[\xi]_{i_d} \dots \dx[\xi]_{i_0}$ is defined recursively as
    \begin{subequations}
        \label{eqn:def-symbol-int}
        \begin{align}
            \int_{t_0}^t \dx[\xi]_{i_d} \dx[\xi]_{i_{d-1}} \dots \dx[\xi]_{i_{0}}
                                        & = \int_{t_0}^t u_{i_d}(\tau) \int_{t_0}^t \dx[\xi]_{i_{d-1}} \dots \dx[\xi]_{i_0} \dx[\tau] \\
            \int_{t_0}^t \dx[\xi]_{i_0} & =  \int_{t_0}^t u_i(\tau) \dx[\tau]
        \end{align}
    \end{subequations}
    with the convention $u_0(t) \equiv 1$.
\end{theorem}
For the system \eqref{eqn:lti-problem}, \eqref{eqn:mycontrol} we calculate the Chen-Fliess expansion with $n= 1, m=2, p=1$ and $h(x) = \operatorname{id}(x)$, $f_i(x)$ as in \eqref{eqn:proof-of-lemma-1-vectorfields} and $u_1(t) = \omega^{\frac 1 2} \sin(\omega t), u_2(t) = \omega^{\frac 1 2} \cos(\omega t)$ for $d = 0, \dots, 3$. We refer to $d$ as the order of the expansion.
\begin{conferenceversion}
    The resulting expressions up to order $d=3$ are given in the extended version \cite{weberebenbauergharesifard2021}.
\end{conferenceversion}
\begin{extendedversion}
    The resulting expressions up to order $d=3$ are given in \Cref{tab:cf-expressions}. The first column displays the identifier as a sequence of indices, where a sequence $s = i_0 \dots i_{d.1} i_d$ corresponds to the indices in \eqref{eqn:cf-input-output-response}. The second column holds the (iterated) Lie-bracket for the same equation, while the third column holds the evaluated integral symbol. The fourth column simplifies the result from the third column by evaluating the Chen-Fliess integral symbol after an integer number of periods of $\omega$, while the last column additionally assumes an initial time $t_0 = 0$.

    We mapped the expressions by their respective indices into tree structures \cref{fig:cf-component-tree-y}, \Cref{fig:cf-component-tree-k} to identify paths which terminate and do not contribute to the output response for higher orders. The trees are to be read in the following way:

    Firstly, there is one tree in \cref{fig:cf-component-tree-y} detailing the contributions to the state $y(t)$, while the second tree in \cref{fig:cf-component-tree-k} contains the contributions to $k(t)$.

    For both trees, every depth level corresponds to one order $d$ of the Chen-Fliess expansion. Every node is structured as
    \begin{center}
        \setlength{\blockheight}{1cm}
        \setlength{\blockdistance}{\blockheight + 1pt}
        \setlength{\blockwidth}{1.3\blockheight -1pt}
        \begin{tikzpicture}[
                grow=right,
                sloped,
                node distance = 0mm and 0mm,
                inner sep = 1pt,
                columnhead/.append style={
                        color=cMBlue,
                    },
                my shape/.style={
                        rectangle split,
                        rectangle split horizontal,
                        rectangle split parts=#1,
                        draw=cLight,
                        minimum height=\blockheight,
                        rectangle split part align=base,
                        every one node part/.style={text width=\blockwidth, align=center},
                        every two node part/.style={text width=\blockwidth, align=center},
                        every three node part/.style={text width=\blockwidth, align=center},
                        every four node part/.style={text width=\blockwidth, align=center},
                        every five node part/.style={text width=\blockwidth, align=center},
                        every six node part/.style={text width=\blockwidth, align=center},
                        fill=white,
                    },
                edge from parent/.append style={
                        draw = cMedium,
                    },
                charting/.append style={
                        child anchor=west,
                        parent anchor=border
                    },
                skipped/.append style={
                        color=cLight,
                    },
                index/.append style={
                        color=black,
                    },
                coefficient/.append style={
                        color = cNormal,
                    },
                powerOfB/.append style={
                        color = cNormal,
                    },
                powerOfY/.append style={
                    },
                powerOfA/.append style={
                    },
                powerOfT/.append style={
                    }]

            \node[my shape=6](ex1){\nodepart[index]{one}$i_d$\nodepart[coefficient, skipped]{two}$c(s)$\nodepart[skipped]{three}$p_b(s)$\nodepart[]{four}$p_y(s)$\nodepart[]{five}$p_r(s)$\nodepart[]{six}$p_T(s)$};
        \end{tikzpicture}
    \end{center}
    where the sequence $s$ needs to be determined by reading the path to the current node from the root and the box contents correspond to terms of the form
    \begin{align*}
        c b^{p_b} y^{p_y}(t_0) (a-bk(t_0))^{p_r} T^{p_t}.
    \end{align*}
    for $T = t- t_0$. Whenever the integral symbol evaluated to zero, we signify this by putting a $\triangle$ in the rightmost box. When the component of the iterated Lie-derivative is zero, we signify this by putting a $\star$ into all but the first to boxes. Once the Lie-derivative is zero for a component, no descendant of this node can be a nonzero value by definition of \eqref{eqn:cf-input-output-response}.

\end{extendedversion}

\section{Proofs}
\label{sec:appendix:proofs}

\begin{proof}[of \texorpdfstring{\Cref{prop:lbs:cartesian}}{Proposition~\ref{prop:lbs:cartesian}}]
    \label{proof:lbs:cartesian}
    We note $u_1(\omega t) = \sin (\omega t)$ and $u_2(\omega t) = \cos(\omega t)$ satisfy \Cref{assumption:ui}.
    Under the conditions of \Cref{prop:lbs:cartesian} and in view of the degrees of freedom used in Lie-bracket approximation, we have $l=2$, ${k_1, k_2 = 1}$, ${p_1, p_2 = \frac 1 2}$, $u_1(k_1 \omega t) = \sin(\omega t)$ and $u_2(k_2 \omega t) = \cos(\omega t)$, which meet \Cref{assumption:ui}, while \Cref{assumption:pi} only applies to the case $p_1 + p_2 > 1$.
    From \eqref{eqn:lti-problem} and \eqref{eqn:mycontrol}, we rearrange terms into vector fields
    \begin{subequations}
        \label{eqn:proof-of-lemma-1-vectorfields}
        \begin{align}
            f_0(t, x) & = \begin{bmatrix}
                              (a-b x_2)\, x_1 \\
                              0
                          \end{bmatrix},
            \\
            f_1(t, x) & = \begin{bmatrix}
                              - b x_1 \\
                              0
                          \end{bmatrix}
            \quad\text{and}               \\
            f_2(t, x) & = \begin{bmatrix}
                              0 \\
                              x_1^2
                          \end{bmatrix},
        \end{align}
    \end{subequations}
    with $x = \begin{bmatrix}
            x_1 & x_2
        \end{bmatrix}^T \coloneqq \begin{bmatrix}
            x & k
        \end{bmatrix}^T$.
    The Lie bracket between vector fields $f_1$ and $f_2$ is calculated as
    \begin{align}
        [f_1, f_2] (t, x) & = \begin{bmatrix}
                                  0 \\ -2 b x_1^2
                              \end{bmatrix}
    \end{align}
    and the coefficient $\gamma_{12}$ in \eqref{eqn:labar:gamma} is independent from $\omega$ by design and can be calculated to be
    \begin{align}
        \gamma_{12} & = \frac{\omega^{\frac 1 2} + \omega^{\frac 1 2}}{2 \pi \omega^{-1}} \!\!\int_0^{\frac{2 \pi}{\omega}} \!\!\! \int_0^\theta \!\!\cos(\omega \theta)  \sin(\omega \tau) \dx[\tau] \dx[\theta] = - \frac 1 2.
    \end{align}
    Hence the Lie-bracket system associated with \eqref{eqn:lti-problem} and \eqref{eqn:mycontrol} takes the very simple form
    \begin{align}
        \label{eqn:tmp:pf1}
        \dot{\bar{x}}(t) & = \begin{bmatrix}
                                 (a- b\bar{x}_2(t)) \,\bar{x}_1(t) \\
                                 b \bar{x}_1^2(t)
                             \end{bmatrix}.
    \end{align}
    Expressing \eqref{eqn:tmp:pf1} in the coordinates $\bar{x}$ and $\bar{k}$ yields the statement of \Cref{prop:lbs:cartesian}.
\end{proof}

\begin{proof}[of \texorpdfstring{\Cref{thm:lbs-system-analysis-cartesian}}{Theorem~\ref{thm:lbs-system-analysis-cartesian}}]
    \label{proof:lbs:properties}
    We introduce the candidate Lyapunov function family
    \vspace{-6pt}
    \begin{align}
        V_p(\bar x) & = \frac 1 2 \bar{x}_1^2 + \frac 1 2 (\bar{x}_2 -c_p)^2, \label{eqn:vc-1}
    \end{align}
    where $c_p = \frac{a+p}{b}$.
    Notice that $V_p$ is positive definite with respect to any point $\bar{x}_p = (0, c_p)$.
    Every function of the family $V_{p}$ is radially unbounded for any fixed $p \in \R$, consequently all sub-level sets
    \vspace{-3pt}
    \begin{align}
        L_p \coloneqq   \left\{ x \in \R^2 : V_p(x) \leq d\right\}
    \end{align}
    are compact. Every function of $V_p$ is at least once continuously differentiable on $\R^2$.
    We take the Lie derivative along solutions of \eqref{eqn:dynamics-lbs-cartesian} to be
    \begin{align}
        \dot{V}_p(\bar{x}(t)) & = \bar{x}_1(t) \,\dot{\bar{x}}_1(t) + (\bar{x}_2(t)-c_p)\,\dot{\bar{x}}_2(t) \\
                              & = (a-c_pb)\,\bar{x}_1^2(t) .\label{eqn:vcdot-1}
    \end{align}
    The set $X_e = \left\{0\right\} \times \R$ is the set of all equilibria of \eqref{eqn:dynamics-lbs-cartesian}, which can be established by calculating the vector field in \eqref{eqn:dynamics-lbs-cartesian}  at any point $x\in X_e$ to be zero, proving \cref{itm:lbs-1}. Moreover, $\bar{x}_p \in X_e$ shows $\bar{x}_p$ is one equilibrium point.
    Note that
    \vspace{-6pt}
    \begin{align}
        \dot{V}_p(\bar{x}(t)) = -p \bar{x}_1^2(t)\leq 0 \label{eqn:vcdot-C0}.
    \end{align}
    Invoking a Lyapunov stability theorem (e.g \cite[Theorem 2.20]{sepulchre2012constructive}) and \eqref{eqn:dynamics-lbs-cartesian} locally Lipschitz, we conclude that $\bar{x}_p$ is a globally stable equilibrium of \eqref{eqn:dynamics-lbs-cartesian}. Since there are other equilibria in any neighbourhood of $\bar{x}_p$, it cannot be attractive. Hence it cannot be asymptotically stable, which concludes \cref{itm:lbs-2}.
    With $d_p = V_p(\bar{x}(t_0))$, define the set
    \begin{align}
        \Omega_p \coloneqq \left\{x \in \R^2 : V_p(x) \leq d_0\right\}.
    \end{align}
    The vector field of \eqref{eqn:dynamics-lbs-cartesian} is time-invariant, thus it  is naturally constant in $t$.
    Moreover, the Jacobian of $f_0 + \gamma_{12}\left[f_1, f_2\right]$ is well-defined for all $\bar{x}\in \R^2$, hence $f_0 + \gamma_{12}\left[f_1, f_2\right]$ is locally Lipschitz in $\bar{x}$ for all $t \geq t_0$ and all $\bar{x} \in \R^2$.
    The set $\Omega_p$ is a sub-level set to $V_p$, hence $\Omega_p$ is compact, i.e. it is closed and bounded. From \eqref{eqn:vcdot-C0} we conclude $V_p(\bar{x}(t)) \leq d_0$. As a result, $\Omega_p$ is shown to be a positively invariant set with respect to \eqref{eqn:dynamics-lbs-cartesian}.
    We conclude with a standard existence and uniqueness theorem \cite[Theorem 3.3]{khalil2002nonlinear}, that for every $x_0\in \R^2$ and for every $t_0 \in \R$, there is a unique solution $\bar{x}(t; t_0, \bar{x}_0)$ to the initial value problem \eqref{eqn:dynamics-lbs-cartesian}, which is defined on $\interval[open right]{t_0}{\infty}$ and is unique. This concludes \cref{itm:lbs-3}

    The set $\Omega_p$ is a disk with radius $d_p$ centered at $\bar{x}_p$. Alternatively, we may write $\Omega_p = \overline{B_{d_0}}(\bar{x}_p)$. Since $\Omega_p$ is a positively invariant set, every solution $\bar{x}(t; t_0, \bar{x}_0)$ lies in $\overline{B_{d_p}}(\bar{x}_p)$ for all $t \geq t_0$, which concludes the statement of \cref{itm:lbs-4}
    Lastly, we will establish the convergence property in \cref{itm:lbs-5}.

    Let $p=1$, then \eqref{eqn:vcdot-1} can be further estimated
    \begin{align}
        \dot{V}_1(\bar{x}(t))  = - \bar{x}_1^2(t) \label{eqn:vcdot-C1}.
    \end{align}
    To continue with LaSalle's theorem, we establish its premises. With $d_1$, define the set
    \begin{align}
        \Omega_1 \coloneqq \left\{x \in \R^2 : V_1(x) \leq d_1\right\}.
    \end{align}
    The set $\Omega_1$ is a sub-level set, therefore $\Omega$ is compact.
    From \eqref{eqn:vcdot-C1}, we conclude $V_{1}(\bar{x}(t)) \leq d$ for all $\bar{x}_0 \in \Omega_1$. As a result, $\Omega_1$ is shown to be a positively invariant set with respect to \eqref{eqn:dynamics-lbs-cartesian}.
    Further, $V_i(\bar{x})$ are continuously differentiable functions and $\dot{V}_{1}(\bar{x}(t)) \leq 0$ on $\Omega_1$ by design and choice of $c_1$. Define according to LaSalle's theorem \cite[Theorem~4.4, p.128]{khalil2002nonlinear}
    \begin{align}
        E
        \coloneqq
        \left\{x \in \Omega_1 : \dot{V}_{1}(x) = 0\right\} = X_e \cap \Omega_1.
    \end{align}
    We conclude with LaSalle's theorem, that every solution $\bar{x}(t; t_0, \bar{x}_0)$ starting in $\Omega_1$ approaches the largest invariant set $M \subseteq E$ as $t \to \infty$.

    This shows convergence of $\bar{x}_1(t)$ to zero as $t$ goes to $\infty$. Moreover, consider only the $\bar{x}_2$-dynamics $\dot{\bar{x}}_2(t) = b\bar{x}_1^2(t)$. Depending on the sign of $b$, the signal $\bar{x}_2(t; t_0, \bar{x}_{2, 0})$ is either monotonically increasing or monotonically decreasing. In either case, it must either converge to a limit or increase respectively decrease to $\pm \infty$. However, since $\Omega_1$ is positively invariant with respect to \eqref{eqn:dynamics-lbs-cartesian} and compact, $\bar{x}_2(t; t_0, \bar{x}_{2, 0})$ must be bounded. Conclusively, $\bar{x}_2(t; t_0, \bar{x}_{2, 0})$ converges towards a limit.
    More can be said about the limit of $\bar{x}_2(t; t_0, \bar{x}_{2,0})$, as we take one more look at $V_0$ and $\Omega_0$.
    We reiterate, $V_i(\bar{x})$ are continuously differentiable functions and $\dot{V}_0(\bar{x}(t)) \leq 0$ in $\Omega_0$ by design and choice of $c_0$. We define once more according to LaSalle's theorem
    \begin{align*}
        E
        \coloneqq
        \left\{x \in \Omega_0 : \dot{V}_{0}(x) = 0\right\} = \left\{x \in \R^2 : V_0(x) = d_0\right\}.
    \end{align*}
    The set $E$ is a circle with radius $d_0$ centered at $\bar{x}_e = \left(0, c_0\right)$. However, the only invariant points on $E$ are
    \begin{align}
        M
        =
        \left\{(0, c_0-d_0), (0, c_0 + d_0)\right\}.
    \end{align}
    Using LaSalle's theorem again, we conclude that any solution $\bar{x}(t; t_0, \bar{x}_0)$ converges to one of two points in $M$.
    Even more can be said about the convergence behaviour. Suppose $b > 0$. Then $\bar{x}_2(t; t_0, \bar{x}_{2,0})$ is monotonically increasing due to \eqref{eqn:dynamics-lbs-cartesian-k} for all $\bar{x}_0 \in \R^2 \setminus X_e$. Consequently, if $b > 0$, then $\lim_{t \to \infty} \bar{x}_2(t; t_0, \bar{x}_{2,0}) = c_0 + d_0$.
    Analogously, assuming $b < 0$ and repeating the previous argument, we conclude that for $ b < 0$, $\lim_{t \to \infty} \bar{x}_2(t; t_0, \bar{x}_{2,0}) = c_0 - d_0$.
\end{proof}

\begin{proof}[of \texorpdfstring{\Cref{prop:polar-coordinate-transformation}}{Proposition~\ref{prop:polar-coordinate-transformation}}]
    \label{proof:polar-coordinates}
    For notational compactness, we apply the following convention: In view of \eqref{eqn:transform-phi}, we denote terms arising from the first case $y \geq 0$ in {\color{cMRed}red} and terms arising from the second case $y < 0$ in {\color{cMBlue}blue}.

    We start by calculating the time derivative of \eqref{eqn:transform-r}
    \begin{align}
        \dot r
         & = \pdv{r}{y} \dot{y} + \pdv{r}{k} \dot{k} + \pdv{r}{t}                                                 \\
         & = \frac{2y \dot{y}}{2\sqrt{y^2 + (k-c_0)^2}} + \frac{(2(k-c_0))\dot{k}}{2\sqrt{y^2 + (k-c_0)^2}} + 0   \\
         & = \frac{1}{r}\left( y\dot{y} + (k-c_0)\dot{k} \right)                                                  \\
         & = \frac{1}{r}\left( y\left[(a - bk)y - b y \cdot \sqrt{\omega} \, u_1(\omega t)\right] \right.\notag   \\
         & \phantom{=\frac 1 r \big(}\left.+ (k-c_0)\left[y^2 \cdot \sqrt{\omega} \, u_2(\omega t)\right] \right)
        \\
        \begin{split}
            &= a r \cos^2 \varphi - br^2 \sin \varphi \cos^2\varphi -b c_0r \cos^2 \varphi \\
            &\phantom{=}- b r  \cos^2 \varphi \cdot \sqrt{\omega} \,u_1(\omega t) \\
            &\phantom{=}+ r^2 \sin\varphi \cos^2 \varphi \cdot \sqrt{\omega} \, u_2(\omega t).
        \end{split}
        \intertext{Using $c_0 = \frac a b$, we calculate further}
        \begin{split}
            &= a r \cos^2 \varphi - br^2 \sin \varphi \cos^2\varphi - a r \cos^2 \varphi \\
            &\phantom{=}- b r  \cos^2 \varphi \cdot \sqrt{\omega} \,u_1(\omega t) \\
            &\phantom{=}+ r^2 \sin\varphi \cos^2 \varphi \cdot \sqrt{\omega} \, u_2(\omega t)
        \end{split}                          \\
        \begin{split}
            &=  - br^2 \sin \varphi \cos^2\varphi \\
            &\phantom{=}- b r  \cos^2 \varphi\cdot \sqrt{\omega}\,  u_1(\omega t) \\
            &\phantom{=}+ r^2 \sin\varphi \cos^2 \varphi \cdot \sqrt{\omega} \, u_2(\omega t)
        \end{split}\label{eqn:r-dot}
    \end{align}
    Likewise, we take the time derivative of \eqref{eqn:transform-phi} to be
    \begin{conferenceversion}%
        \begin{align}
            \dot \varphi
            =   & {} \pdv{\varphi}{y}\dot{y} + \pdv{\varphi}{k}\dot{k} + \pdv{\varphi}{t}         \\
            \begin{split}
                =
                &\hphantom{+} \pdv{}{y}\left({\color{cMBlue}\pi -} \arcsin \frac {k-c_0} r\right) \dot{y} \\
                &+ \pdv{}{k}\left({\color{cMBlue}\pi -} \arcsin \frac {k-c_0} r\right) \dot{k} \\
                &+ \,0
            \end{split}
            \\
            =   & A_0 \cdot A_1, \label{eqn:a0-a1-split}
            \intertext{where}
            A_0 & = \cpm \frac{1}{\sqrt{1-\left(\frac {k-c_0} r\right)^2}}\quad \text{and} \notag \\
            A_1 & = \frac{ - \frac{1}{2r}(2y)\cdot (k-c_0)}{r^2} \dot{y}
            + \frac{r- \frac{1}{2r}2(k-c_0)^2}{r^2} \dot{k}. \notag
            \intertext{After lengthy calculations, detailed in the extended version, we have}
            A_0
                & =
            \frac{1}{\cos \varphi} \quad \text{and}\label{eqn:a0-final},
            \\
            A_1  \begin{split}
                     &= \phantom{+} b  r \sin^2 \varphi \cos^2 \varphi                             \\
                     &\phantom{=} +  b  \sin \varphi \cos^2 \varphi     \cdot \sqrt{\omega} \, u_1(\omega t)  \\
                     &\phantom{=} +  r \cos^4 \varphi                  \cdot \sqrt{\omega} \, u_2(\omega t) .
                 \end{split}
            %
            %
            \label{eqn:a1-final}
        \end{align}
        Resubstituting \eqref{eqn:a0-final} and \eqref{eqn:a1-final} into \eqref{eqn:a0-a1-split}, we are left with the transformed dynamics
        \begin{align}
            \dot{\varphi} =
            \begin{split}
                & \phantom{-} b  r \sin^2 \varphi \cos \varphi                                    \\
                &+b  \sin \varphi \cos \varphi     \cdot\sqrt{\omega} \, u_1(\omega t)  \\
                & + r \cos^3 \varphi                \cdot \sqrt{\omega} \, u_2(\omega t).
            \end{split}
            \label{eqn:phi-dot}
        \end{align}
    \end{conferenceversion}%
    \begin{extendedversion}

        \begin{align}
            \dot \varphi
            =   & {} \pdv{\varphi}{y}\dot{y} + \pdv{\varphi}{k}\dot{k} + \pdv{\varphi}{t}                                                         \\
            \begin{split}
                =
                &\hphantom{+} \pdv{}{y}\left({\color{cMBlue}\pi -} \arcsin \frac {k-c_0} r\right) \dot{y} \\
                &+ \pdv{}{k}\left({\color{cMBlue}\pi -} \arcsin \frac {k-c_0} r\right) \dot{k} \\
                &+ \,0
            \end{split}
            \\
            \begin{split}
                = & \cpm \frac{1}{\sqrt{1-\left(\frac {k-c_0} r\right)^2}}
                \cdot
                \frac{0 - \frac{\partial r}{\partial y}\cdot (k-c_0)}{r^2} \dot{y}
                \\
                &
                \cpm \frac{1}{\sqrt{1-\left(\frac {k-c_0} r\right)^2}}
                \cdot
                \frac{ 1\cdot r - \frac{\partial r}{\partial k}\cdot(k-c_0)}{r^2} \dot{k}
            \end{split}
            \\
            =   & A_0 \cdot A_1, \label{eqn:a0-a1-split}
            \intertext{where}
            A_0 & = \cpm \frac{1}{\sqrt{1-\left(\frac {k-c_0} r\right)^2}}\quad \text{and} \notag                                                 \\
            A_1 & = \frac{ - \frac{1}{2r}(2y)\cdot (k-c_0)}{r^2} \dot{y}
            + \frac{r- \frac{1}{2r}2(k-c_0)^2}{r^2} \dot{k}. \notag
            \intertext{Calculating further, we have}
            A_0
                & =
            \cpm \frac{1}{\sqrt{1-\left(\frac {(r \sin \varphi + c_0)-c_0} r\right)^2}}                                                           \\
                & =
            \cpm \frac{1}{\sqrt{1-\left(\frac {r \sin \varphi } r\right)^2}}                                                                      \\
                & =
            \cpm \frac{1}{\sqrt{1- \sin^2 \varphi }}                                                                                              \\
                & =
            \cpm \frac{1}{|\cos \varphi|} \label{eqn:a0-before-case-merge}                                                                        \\
                & =
            \frac{1}{\cos \varphi}\label{eqn:a0-final},
            \intertext{where from \eqref{eqn:a0-before-case-merge} to \eqref{eqn:a0-final} we used the fact that $y \geq 0$ implies $\operatorname{sign}(\cos \varphi) = 1$ and $y < 0$ implies $\operatorname{sign}(\cos \varphi) = -1$. We calculate the second factor to be}
            A_1 & =
            \frac{1}{r^2}\cdot \Big(
            - \frac{1}{r}y (k-c_0) \dot{y}
            + (r- \frac{1}{r}(k-c_0)^2)\dot{k}
            \Big)
            \\
                & = \frac{1}{r^3}\cdot \Big(
            - y (k-c_0) \dot{y}
            + (r^2- (k-c_0)^2)\dot{k}
            \Big)
            \\
                & =\frac{1}{r^3} \cdot \left(
            - y (k-c_0) \dot{y}
            + y^2\dot{k}
            \right)
            \\
            %
                & \begin{aligned}[b]
                      {}=  \frac{1}{r^3} \cdot \Big(
                       & - y(k-c_0)  \mathrlap{\left[(a - bk)y - b y \cdot \sqrt{\omega} \, u_1(\omega t) \right]} \\
                       & + y^2       \left[y^2 \cdot \sqrt{\omega} \, u_2(\omega t)\right]\quad\Big)
                  \end{aligned}
            \\
                & \begin{aligned}[b]%
                      {}=\frac{1}{r^3} \cdot \Big( 
                       & -  (k-c_0)(a-bk)y^2                                                  \\
                       & +b(k-c_0)y^2         \cdot \sqrt{\omega} \, u_1(\omega t)            \\
                       & + y^4                \cdot \sqrt{\omega} \, u_2(\omega t) \quad\Big)
                  \end{aligned}
            \\
                & \begin{aligned}[b]
                      {}=  \frac{1}{r^3} \cdot \Big(
                       & -  ((r\sin \varphi + c_0)-c_0)\cdot \dots                                                                                    \\
                       & \qquad  (a-b(r \sin \varphi + c_0))
                      (r \cos \varphi)^2                                                                                                              \\
                       & +b((r \sin \varphi + c_0)-c_0)(r \cos \varphi)^2                            \cdot \dots                                      \\
                       & \qquad \sqrt{\omega} \, u_1(\omega t)                                                                                        \\
                       & + (r \cos \varphi)^4                                                        \cdot \sqrt{\omega} \, u_2(\omega t) \quad \Big)
                  \end{aligned}
            \\
                & \begin{aligned}[b]
                      {}=  \frac{1}{r^3}  \cdot \Big(
                       & -  (r\sin \varphi )(a-b(r \sin \varphi + c_0))\mathrlap{(r \cos \varphi)^2}                                       \\
                       & +b(r \sin \varphi )(r \cos \varphi)^2                             \cdot \sqrt{\omega} \, u_1(\omega t)            \\
                       & + (r \cos \varphi)^4                                              \cdot \sqrt{\omega} \, u_2(\omega t)\quad \Big)
                  \end{aligned}
            \\
                & \begin{aligned}[b]
                      {}=  \frac{1}{r^3} \cdot \Big(
                       & -  (a-b(r \sin \varphi + \frac a b))r^3 \mathrlap{\sin \varphi \cos^2 \varphi  }                                      \\
                       & +br^3  \sin \varphi \cos^2 \varphi                                   \cdot \sqrt{\omega} \, u_1(\omega t)             \\
                       & + r^4 \cos^4 \varphi                                                 \cdot \sqrt{\omega} \, u_2(\omega t) \quad \Big)
                  \end{aligned}
            \\
                & \begin{aligned}[b]
                      {}=\frac{1}{r^3} \cdot \Big(
                       & -  (a-br \sin \varphi -  a )r^3 \mathrlap{\sin \varphi \cos^2 \varphi}                                        \\
                       & +br^3  \sin \varphi \cos^2 \varphi                           \cdot \sqrt{\omega} \, u_1(\omega t)             \\
                       & + r^4 \cos^4 \varphi                                         \cdot \sqrt{\omega} \, u_2(\omega t) \quad \Big)
                  \end{aligned}
            \\
                & \begin{aligned}[b]
                      {}= \frac{1}{r^3} \cdot \Big(
                       & \phantom{+} b  r^4 \sin^2 \varphi \cos^2 \varphi                                                  \\
                       & + br^3  \sin \varphi \cos^2 \varphi               \cdot \sqrt{\omega} \, u_1(\omega t)            \\
                       & + r^4 \cos^4 \varphi                             \cdot \sqrt{\omega} \, u_2(\omega t) \quad \Big)
                  \end{aligned}
            \\
                & \begin{aligned}[b]
                      {}=
                       & \hphantom{+}  b  r \sin^2 \varphi \cos^2 \varphi                        \\
                       &
                      +  b  \sin \varphi \cos^2 \varphi     \cdot \sqrt{\omega} \, u_1(\omega t) \\
                       &
                      +  r \cos^4 \varphi                  \cdot \sqrt{\omega} \, u_2(\omega t)
                  \end{aligned}                                                      \\
            %
            %
            %
                & \begin{aligned}[b]
                      {}=
                       & \hphantom{+} b  r \sin^2 \varphi \cos^2 \varphi                             \\
                       & +  b  \sin \varphi \cos^2 \varphi     \cdot \sqrt{\omega} \, u_1(\omega t)  \\
                       & +  r \cos^4 \varphi                  \cdot \sqrt{\omega} \, u_2(\omega t) .
                  \end{aligned}
            \label{eqn:a1-final}
        \end{align}
        Resubstituting \eqref{eqn:a0-final} and \eqref{eqn:a1-final} into \eqref{eqn:a0-a1-split}, we are left with the transformed dynamics
        \begin{align}
            \dot{\varphi} =
            \begin{split}
                & \phantom{-} b  r \sin^2 \varphi \cos \varphi                                    \\
                &+b  \sin \varphi \cos \varphi     \cdot\sqrt{\omega} \, u_1(\omega t)  \\
                & + r \cos^3 \varphi                \cdot \sqrt{\omega} \, u_2(\omega t).
            \end{split}
            \label{eqn:phi-dot}
        \end{align}
    \end{extendedversion}

    \noindent With the same parameter as in the proof of \Cref{prop:lbs:cartesian}, we rewrite \eqref{eqn:r-dot} and \eqref{eqn:phi-dot} in the form of \eqref{eqn:labar:system} using the vector fields
    \begin{align}
        f_0 & = 
 $.
    The Lie-bracket between $f_1$ and $f_2$ can be calculated as
    \begin{align}
        [f_1, f_2]
         & =
        -2 b x_1\begin{aligned}[t]
                     & %
                \end{aligned}
    \end{align}
    and  the Lie-bracket system is
    \begin{align}
        %
.
    \end{align}
    Renaming $\bar{r}(t) = \bar{x}_1(t)$ and $\bar{\varphi}(t) = \bar{x}_2(t)$ completes the proof.
\end{proof}

\begin{extendedversion}
    \onecolumn%
    \section{Calculations}
    \input{extended/cf_calculations.tex}%
    \setlength{\topmargin}{-2cm}
    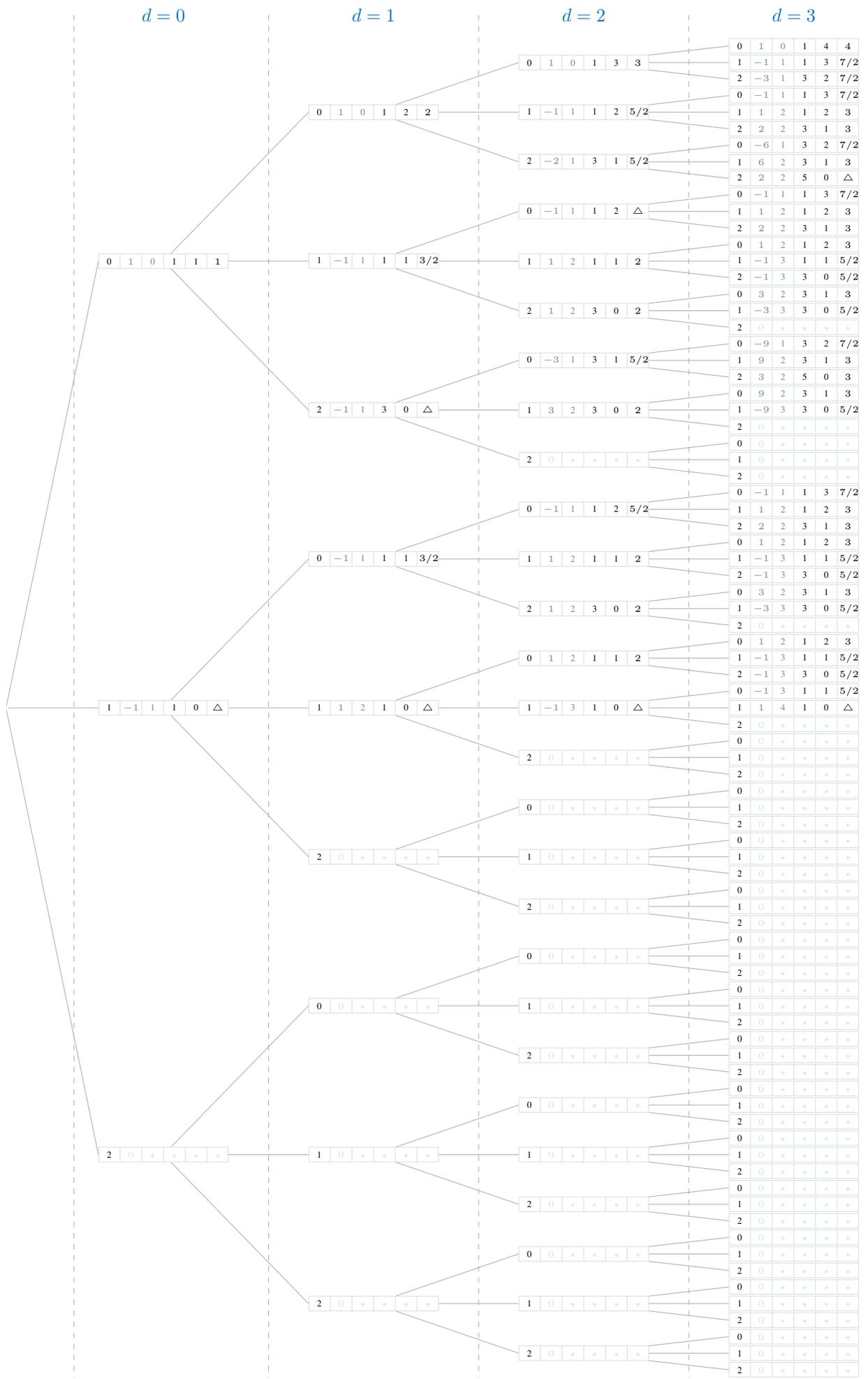
\begin{figure}[H]%
        \tikzsetnextfilename{figure_cf_component1}%
        \label{fig:cf-component-tree-y}
        \centering%
        \input{figures/tree_y.tikz}%
        \caption{Tree chart of the expressions involved in the Chen-Fliess expansion for expressions contributing to $y(t)$.}
    \end{figure}%

    \begin{figure}[H]%
        \tikzsetnextfilename{figure_cf_component2}%
        \label{fig:cf-component-tree-k}
        \centering%
        \input{figures/tree_k.tikz}%
        \caption{Tree chart of the expressions involved in the Chen-Fliess expansion for expressions contributing to $k(t)$.}
    \end{figure}%

    \twocolumn
\end{extendedversion}

\end{document}

%% file: helper_setup.tex
\definecolor{cNormal}{HTML}{7d8391} 
\definecolor{cLight}{HTML}{d8e0e3} 
\colorlet{cMedium}{cLight!50!cNormal} 

\definecolor{cMBlue}{rgb}{0, 0.4470, 0.7410}
\definecolor{cMOrange}{rgb}{0.8500, 0.3250, 0.0980}
\definecolor{cMYellow}{rgb}{0.9290, 0.6940, 0.1250}
\definecolor{cMViolet}{rgb}{0.4940, 0.1840, 0.5560}
\definecolor{cMGreen}{rgb}{0.4660, 0.6740, 0.1880}
\definecolor{cMBabyblue}{rgb}{0.3010, 0.7450, 0.9330}
\definecolor{cMRed}{rgb}{0.6350, 0.0780, 0.1840}

\colorlet{cdemoversion}{cMYellow}
\colorlet{cextendedversion}{cMViolet}
\colorlet{cconferenceversion}{cMGreen}

\intervalconfig{
    left open fence ={(},
    right open fence = {)}
}

\hypersetup{linkbordercolor=green}

\bookmarksetup{open, numbered, startatroot}

\makeatletter
\newcommand{\thlabel}[1]{
    \thmt@thmname\space
    \ifx\@currentlabelname\@empty
    \else
        \space\@currentlabelname%
    \fi
}
\newcommand*{\theorembookmark}{%
    \bookmark[
        dest=\@currentHref,
        rellevel=1,
        keeplevel,
    ]{\thlabel{\thmt@envname}}
}
\makeatother

\allowdisplaybreaks

\makeatletter
\if@confmode%
    \def\subsection{\@startsection{subsection}{4}{\z@}{0ex plus 0.1ex minus 0.1ex}%
        {0ex}{\normalfont\normalsize\itshape}}%
\else 
    \def\subsection{\@startsection{subsection}{2}{\z@}{3.5ex plus 1.5ex minus 1.5ex}%
        {0.7ex plus .5ex minus 0ex}{\normalfont\normalsize\itshape}}%
\fi

\def\@IEEEprocessthesectionargument#1{%
    \@ifmtarg{#1}{%
        \@IEEEappendixsavesection*{Appendix \thesectiondis}%
        \addcontentsline{toc}{section}{Appendix \thesection}}{%
        \@IEEEappendixsavesection*{Appendix \thesectiondis: #1}%
        \addcontentsline{toc}{section}{Appendix \thesection: #1}}}
\makeatother

\makeatletter
\declaretheoremstyle[
    spaceabove=6pt,
    spacebelow=6pt,
    bodyfont=\normalfont\itshape,
    numbered=yes,
    qed=\qedsymbol,
    headformat=\NAME \NOTE ,
    notebraces={}{},
    notefont=\normalfont\bfseries,
]{proofstyle}
\makeatother

\declaretheoremstyle[
    spaceabove=6pt,
    spacebelow=6pt,
    bodyfont=\normalfont\itshape,
    headfont=\normalfont\bfseries,
    headpunct={},
]{theoremstyle}

\makeatletter
\let\c@theorem\relax
\makeatother

\declaretheorem[name=Theorem, postheadhook=\theorembookmark, style=theoremstyle]{theorem}
\declaretheorem[name=Assumption, postheadhook=\theorembookmark, style=theoremstyle]{assumption}
\declaretheorem[name=Proposition, postheadhook=\theorembookmark, style=theoremstyle]{proposition}

\makeatletter
\let\c@lemma\relax
\makeatother

\declaretheorem[name=Lemma, postheadhook=\theorembookmark, style=theoremstyle]{lemma}

\declaretheorem[name = Proof, postheadhook=\theorembookmark, style=proofstyle]{proof}

\newlist{assenum}{enumerate}{1} 
\setlist[assenum]{label=\normalfont\arabic*., ref={\theassumption-\arabic*}}
\crefalias{assenumi}{assumption}
\crefname{assumption}{assumption}{assumptions}
\Crefname{assumption}{Assumption}{Assumptions}
\crefname{statement}{statement}{statements}
\Crefname{statement}{Statement}{Statements}
\crefname{paragraph}{paragraph}{paragraphs}
\Crefname{paragraph}{Paragraph}{Paragraphs}

\newcommand{\R}{\mathbb{R}}%
\newcommand{\Q}{\mathbb{Q}}%
\newcommand{\Rnot}[1]{\mathbb{R}\!\setminus\!\{#1\}\xspace}%

\newcommand{\dx}[1][]{%
\ifthenelse{ \equal{#1}{} }
{\ensuremath{\operatorname{d}\!{x}}}
{\ensuremath{\operatorname{d}\!{#1}}}
}

\newcommand{\LieDeriv}[3][]{%
    \ifthenelse{ \equal{#1}{} }
    {\ensuremath{\mathcal{L}_{#2}{#3}}}
    {\ensuremath{\mathcal{L}_{#2}{#3}(#1)}}
}

\newcommand{\LieB}[3][]{%
    \ifthenelse{ \equal{#1}{} }
    {\ensuremath{\[{#2}, {#3}\]}}
    {\ensuremath{[{#2}, {#3}](#1)}}
}

\def\cpm{\mathbin{\ThisStyle{\ensurestackMath{\abovebaseline[-\dimexpr1pt+2.4\LMpt]{%
                    \stackunder[-\dimexpr1pt+2.5\LMpt]{\color{cMRed}\SavedStyle+}{%
                        \color{cMBlue}\SavedStyle-}}}}}}

\newlength{\blockheight}
\newlength{\blockdistance}
\newlength{\blockwidth}

\pgfplotsset{compat=newest}

\pgfdeclarelayer{background}
\pgfsetlayers{background,main}

\pgfplotsset{
    colormap={cm1}{
            color=(cMBlue)
            color=(cMRed)
        },
    inner axis line style = {
            draw=cNormal,
        },
    outer axis line style = {
            draw=cNormal,
        },
    every axis/.append style={
            grid=both,
            grid style={
                    line width=.1pt,
                    draw=cLight,
                },
            major grid style={
                    line width=.2pt,
                    draw=cMedium,
                },
            minor tick num=5,
            ticklabel style={
                    font=\tiny,
                    fill=white,
                    inner sep= 1pt,
                },
            extra y tick style={%
                    grid=major,
                },
            scale only axis,
        },
    lbsplot/.append style={
            densely dashed,
            draw opacity = 0.75,
        },
    esplot/.append style = {
        },
    icplot1/.append style = {
            draw = cMBlue,
        },
    icplot2/.append style = {
            draw = cMRed,
        },
    icplot3/.append style = {
            draw = cMBabyblue,
        },
    icplot4/.append style = {
            draw = cMOrange,
        },
    icplot5/.append style = {
            draw = cMRed,
            mark size=1pt,
            mark = x,
            mark repeat={8},
        },
    icplot6/.append style = {
            draw = cMGreen,
            mark size=1pt,
        },
    icplot7/.append style = {
            draw = cNormal,
            mark = o,
            mark repeat={8},
            mark size=1pt,
        },
    icplot8/.append style = {
            draw = cMOrange,
            mark=x,
            mark repeat={8},
            mark size=1pt,
        },
    ourplot/.append style = {
            color = cMBlue,
            mark size = 1pt,
            mark repeat = {8}
        },
    otherplot/.append style = {
            color = cMRed,
            mark size = 1pt,
            mark repeat = {8}
        },
    var1/.append style= {
        },
    var2/.append style= {
            scatter=true,
            mark = +,
            visualization depends on = {min(abs(x/2), 1pt) \as \perpointmarksize},
            scatter/@pre marker code/.style={/tikz/mark size=\perpointmarksize},
            scatter/@post marker code/.style={},
        },
    var3/.append style= {
            scatter=true,
            mark = square,
            visualization depends on = {min(abs(x/2), 0.5pt) \as \perpointmarksize},
            scatter/@pre marker code/.style={/tikz/mark size=\perpointmarksize},
            scatter/@post marker code/.style={},
        },
    var4/.append style= {
            scatter=true,
            mark = o,
            visualization depends on = {min(abs(x/2), 0.5pt) \as \perpointmarksize},
            scatter/@pre marker code/.style={/tikz/mark size=\perpointmarksize},
            scatter/@post marker code/.style={},
        },
    varlight/.append style = {
            draw = cMBabyblue,
        },
    phasediagram/.style={
            axis lines = middle,
            axis equal image,
            scale only axis,
            axis line style={-Latex},
            xlabel style={at={(ticklabel* cs:1)},anchor=south},
            ylabel style={at={(ticklabel* cs:1)},anchor=west },
            xlabel={$y$},
            ylabel={$k$},
            view     = {0}{90},
        },
    every axis y label/.append style ={
            rotate=-90,
        },
    every tick label/.append style={
            font = \tiny,
        },
    every x tick/.append style={cNormal},
    every y tick/.append style={cNormal},
    StyleVectorfield/.style={%
            color=cLight,
        },
    groupplotstyle2by1/.style={
            group style= {
                    group size = 2 by 1,
                    xlabels at=edge bottom,
                    ylabels at=edge left,
                    vertical sep =0.5cm,
                    horizontal sep =0.75cm,
                },
            xlabel={$t$},
            clip = false,
            xmin=0,
            xmax = 1,
            height=2.5cm,
            width=\textwidth/2-0.75cm/2,
        },
    groupplotstyle1by2/.style={
            group style= {
                    group size = 1 by 2,
                    xlabels at=edge bottom,
                    ylabels at=edge left,
                    vertical sep =0.5cm,
                },
            xlabel={$t$},
            xmin=0,
            xmax = 0.5,
            height=2.5cm,
            width=\linewidth-0.5cm
        },
    groupplotstyle2by1small/.style={
            group style= {
                    group size = 2 by 1,
                    xlabels at=edge bottom,
                    ylabels at=edge left,
                    horizontal sep =0.5cm,
                },
            xmin=0,
            xmax = 0.5,
            height=2.5cm,
            width=\linewidth/2-0.5cm/2,
        }
}

\excludeversion{conferenceversion}
\includeversion{extendedversion}
\excludeversion{demoversion}

\makeatletter
\renewcommand\beginmarkversion{\color{c\@currenvir}}%
\makeatother


%% file: figures/figure_orbits.tikz
\begin{tikzpicture}[
        trim axis left,
        trim axis right,
        background rectangle/.style={fill=white},
        show background rectangle,
        inner frame sep = .25cm,
        every node/.append style={overlay},
    ]%
    \pgfplotstableread[%
        skip first n = 8, 
        col sep =comma,%
        header = true,%
    ]{simdata/works_quad1_3.csv}{\orbitA}%
    \pgfplotstableread[%
        skip first n = 8, 
        col sep =comma,%
        header = true,%
    ]{simdata/works_quad2_1.csv}{\orbitB}%
    \pgfplotstableread[%
        skip first n = 8, 
        col sep =comma,%
        header = true,%
    ]{simdata/works_quad2_2.csv}{\orbitC}%
    \pgfplotstableread[%
        skip first n = 8, 
        col sep =comma,%
        header = true,%
    ]{simdata/works_quad2_3.csv}{\orbitD}%
    %
    %
    %
    \pgfmathsetmacro\PGFa{10}%
    \pgfmathsetmacro\PGFb{-2}%
    \pgfmathsetmacro\PGFczero{(\PGFa / \PGFb)}%
    \pgfmathsetmacro\PGFcone{(\PGFa / \PGFb + 1/\PGFb)}%
    %
    \pgfmathsetmacro\PGFYL{-9}%
    \pgfmathsetmacro\PGFYU{-1}%
    \pgfmathsetmacro\PGFKL{-14}%
    \pgfmathsetmacro\PGFKU{4}%
    %
    \pgfmathtruncatemacro\IY{(\PGFYU-\PGFYL+1)}
    \pgfmathtruncatemacro\IK{(\PGFKU-\PGFKL+1)}

    \pgfplotstableset{
        create on use/y/.style={
                create col/expr={
                        (\PGFYL+\pgfplotstablerow)
                    }
            },
        create on use/k/.style={
                create col/expr={
                        \I
                    }
            }
    }

    \pgfplotsforeachungrouped \i in {(\PGFKL),(\PGFKL+1),...,(\PGFKU)} {%
            \pgfmathsetmacro\I{\i}
            \pgfplotstablenew[columns={y, k}]{\IY}\newpart
            \pgfplotstablevertcat{\grid}{\newpart} 
        }%


    \begin{axis} [
            xmin = -9.99, xmax = 9.99,
            ymin = -12.99, ymax = 1.99,
            phasediagram,
            width=\linewidth-0.5cm,
            height=\linewidth-0.5cm
        ]

        \pgfmathdeclarefunction{dy}{2}{\pgfmathparse{ (\PGFa-\PGFb*#2)*#1}}
        \pgfmathdeclarefunction{dk}{2}{\pgfmathparse{ \PGFb*#1^2}}
        \pgfmathdeclarefunction{n2}{2}{\pgfmathparse{ sqrt((dy(#1, #2))^2+(dk(#1, #2))^2) }}

        \pgfkeys{/pgf/fpu=true}
        \pgfmathsetmacro\PGFmax{(n2(\PGFYL,\PGFKL))}
        \pgfkeys{/pgf/fpu=false}

        \addplot[%
        StyleVectorfield,
        point meta={n2(x,y)},
        quiver={
        u={dy(x, y)/\PGFmax},
        v={dk(x,y)/\PGFmax},
        scale arrows=1,
        every arrow/.append style={%
        line width={max(0.1, \pgfplotspointmetatransformed/1000)},
        -{Stealth[length=0pt 3, width=0pt 3]},
        },
        },
        ] table [x expr=(\thisrow{y}), y expr=(\thisrow{k})] {\grid};

        \addplot[%
        StyleVectorfield,
        point meta={n2(x,y)},
        quiver={
        u={dy(x, y)/\PGFmax},
        v={dk(x,y)/\PGFmax},
        scale arrows=1,
        every arrow/.append style={%
        line width={max(0.1, \pgfplotspointmetatransformed/1000)},
        -{Stealth[length=0pt 3, width=0pt 3]},
        },
        },
        ] table [x expr=(-\thisrow{y}), y expr=(\thisrow{k})] {\grid};

        \coordinate[
        ] (c0) at (axis cs:0,\PGFczero,0) ;
        \fill (c0) circle (1pt);



        \addplot [lbsplot, ourplot] table[
                x = lbs.y,
                y = lbs.k
            ] {\orbitA};
        \addplot [esplot, var2, ourplot] table[
                x = es.y,
                y = es.k
            ] {\orbitA};

        \addplot [lbsplot, ourplot, var1, varlight ] table[
                x = lbs.y,
                y = lbs.k
            ] {\orbitB};
        \addplot [esplot, ourplot, varlight] table[
                x = es.y,
                y = es.k
            ] {\orbitB};

        \addplot [lbsplot, ourplot] table[
                x = lbs.y,
                y = lbs.k
            ] {\orbitC};
        \addplot [esplot, ourplot, var3] table[
                x = es.y,
                y = es.k
            ] {\orbitC};

        \addplot [lbsplot, ourplot] table[
                x = lbs.y,
                y = lbs.k
            ] {\orbitD};
        \addplot [esplot,
            ourplot,
            var4] table[
                x = es.y,
                y = es.k
            ] {\orbitD};
    \end{axis}
\end{tikzpicture}%

%% file: figures/figure_orbits_swapped.tikz
\begin{tikzpicture}[
        background rectangle/.style={fill=white},
        show background rectangle,
        inner frame sep = .25cm,
        every node/.append style={overlay},
    ]

    \pgfplotstableread[
        skip first n = 8, 
        col sep =comma,
        header = true,
    ]{simdata/conference_swapped_system.csv}{\orbit}

    \pgfmathsetmacro\PGFa{10}
    \pgfmathsetmacro\PGFb{-2}
    \pgfmathsetmacro\PGFczero{(\PGFa / \PGFb)}
    \pgfmathsetmacro\PGFcone{(\PGFa / \PGFb + 1/\PGFb)}

    \pgfmathsetmacro\PGFYL{-11}
    \pgfmathsetmacro\PGFYU{-1}
    \pgfmathsetmacro\PGFKL{-14}
    \pgfmathsetmacro\PGFKU{4}

    \pgfmathtruncatemacro\IY{(\PGFYU-\PGFYL+1)}
    \pgfmathtruncatemacro\IK{(\PGFKU-\PGFKL+1)}

    \pgfplotstableset{
        create on use/y/.style={
                create col/expr={
                        (\PGFYL+\pgfplotstablerow)
                    }
            },
        create on use/k/.style={
                create col/expr={
                        \I
                    }
            }
    }

    \pgfplotsforeachungrouped \i in {(\PGFKL),(\PGFKL+1),...,(\PGFKU)} {%
            \pgfmathsetmacro\I{\i}
            \pgfplotstablenew[columns={y, k}]{\IY}\newpart
            \pgfplotstablevertcat{\grid}{\newpart} 
        }%


    \begin{axis} [
            xmin = -7.99, xmax = 11.99,
            ymin = -12.99, ymax = 2.99,
            phasediagram,
            width=\linewidth-0.5cm,
            height=\linewidth-0.5cm
        ]

        \pgfmathdeclarefunction{dy}{2}{\pgfmathparse{ (\PGFa-\PGFb*#2)*#1}}
        \pgfmathdeclarefunction{dk}{2}{\pgfmathparse{ \PGFb*#1^2}}
        \pgfmathdeclarefunction{n2}{2}{\pgfmathparse{ sqrt((dy(#1, #2))^2+(dk(#1, #2))^2) }}

        \pgfkeys{/pgf/fpu=true}
        \pgfmathsetmacro\PGFmax{(n2(\PGFYL,\PGFKL))}
        \pgfkeys{/pgf/fpu=false}

        \addplot[%
        StyleVectorfield,
        point meta={n2(x,y)},
        quiver={
        u={dy(x, y)/\PGFmax},
        v={dk(x,y)/\PGFmax},
        scale arrows=1,
        every arrow/.append style={%
        line width={max(0.1, \pgfplotspointmetatransformed/1000)},
        -{Stealth[length=0pt 3, width=0pt 3]},
        },
        },
        ] table [x expr=(\thisrow{y}), y expr=(\thisrow{k})] {\grid};

        \addplot[%
        StyleVectorfield,
        point meta={n2(x,y)},
        quiver={
        u={dy(x, y)/\PGFmax},
        v={dk(x,y)/\PGFmax},
        scale arrows=1,
        every arrow/.append style={%
        line width={max(0.1, \pgfplotspointmetatransformed/1000)},
        -{Stealth[length=0pt 3, width=0pt 3]},
        },
        },
        ] table [x expr=(-\thisrow{y}), y expr=(\thisrow{k})] {\grid};

        \coordinate[
        ] (c0) at (axis cs:0,\PGFczero,0) ;
        \fill (c0) circle (1pt);



        \addplot [lbsplot, ourplot] table[
                x = lbs.y,
                y = lbs.k
            ] {\orbit};\label{plot:ex-swapped:lbs}
        \addplot [esplot, ourplot] table[
                x = es.y,
                y = es.k
            ] {\orbit};\label{plot:ex-swapped:es}

        \addplot [esplot, otherplot] table[
                x = es_swap.y,
                y = es_swap.k
            ] {\orbit};\label{plot:ex-swapped:es-swapped}

    \end{axis}

\end{tikzpicture}

%% file: figures/figure_orbits_cf.tikz
\begin{tikzpicture}[%
        background rectangle/.style={fill=white},
        show background rectangle,
        inner frame sep = .25cm,
        every node/.append style={overlay},
    ]%

    \pgfplotstableread[
        skip first n = 8, 
        col sep =comma,
        header = true,
    ]{simdata/conference_cf.csv}{\orbit}

    \pgfmathsetmacro\PGFa{10}
    \pgfmathsetmacro\PGFb{-2}
    \pgfmathsetmacro\PGFczero{(\PGFa / \PGFb)}
    \pgfmathsetmacro\PGFcone{(\PGFa / \PGFb + 1/\PGFb)}

    \pgfmathsetmacro\PGFYL{-9}
    \pgfmathsetmacro\PGFYU{-1}
    \pgfmathsetmacro\PGFKL{-14}
    \pgfmathsetmacro\PGFKU{4}

    \pgfmathtruncatemacro\IY{(\PGFYU-\PGFYL+1)}
    \pgfmathtruncatemacro\IK{(\PGFKU-\PGFKL+1)}

    \pgfplotstableset{
        create on use/y/.style={
                create col/expr={
                        (\PGFYL+\pgfplotstablerow)
                    }
            },
        create on use/k/.style={
                create col/expr={
                        \I
                    }
            }
    }

    \pgfplotsforeachungrouped \i in {(\PGFKL),(\PGFKL+1),...,(\PGFKU)} {%
            \pgfmathsetmacro\I{\i}
            \pgfplotstablenew[columns={y, k}]{\IY}\newpart
            \pgfplotstablevertcat{\grid}{\newpart} 
        }%


    \begin{axis} [
            xmin = -9.99, xmax = 9.99,
            ymin = -12.99, ymax = 2.99,
            phasediagram,
            width=\linewidth-0.5cm,
            height=\linewidth-0.5cm
        ]

        \pgfmathdeclarefunction{dy}{2}{\pgfmathparse{ (\PGFa-\PGFb*#2)*#1}}
        \pgfmathdeclarefunction{dk}{2}{\pgfmathparse{ \PGFb*#1^2}}
        \pgfmathdeclarefunction{n2}{2}{\pgfmathparse{ sqrt((dy(#1, #2))^2+(dk(#1, #2))^2) }}

        \pgfkeys{/pgf/fpu=true}
        \pgfmathsetmacro\PGFmax{(n2(\PGFYL,\PGFKL))}
        \pgfkeys{/pgf/fpu=false}

        \addplot[%
        StyleVectorfield,
        point meta={n2(x,y)},
        quiver={
        u={dy(x, y)/\PGFmax},
        v={dk(x,y)/\PGFmax},
        scale arrows=1,
        every arrow/.append style={%
        line width={max(0.1, \pgfplotspointmetatransformed/1000)},
        -{Stealth[length=0pt 3, width=0pt 3]},
        },
        },
        ] table [x expr=(\thisrow{y}), y expr=(\thisrow{k})] {\grid};

        \addplot[%
        StyleVectorfield,
        point meta={n2(x,y)},
        quiver={
        u={dy(x, y)/\PGFmax},
        v={dk(x,y)/\PGFmax},
        scale arrows=1,
        every arrow/.append style={%
        line width={max(0.1, \pgfplotspointmetatransformed/1000)},
        -{Stealth[length=0pt 3, width=0pt 3]},
        },
        },
        ] table [x expr=(-\thisrow{y}), y expr=(\thisrow{k})] {\grid};

        \coordinate[
        ] (c0) at (axis cs:0,\PGFczero,0) ;
        \fill (c0) circle (1pt);



        \addplot [lbsplot, ourplot] table[
                x = lbs.y,
                y = lbs.k
            ] {\orbit};\label{plot:ex3:lbs}
        \addplot [esplot, ourplot, var4] table[
                x = es.y,
                y = es.k
            ] {\orbit};\label{plot:ex3:es}

        \addplot [esplot, otherplot] table[
                x = cf0.y,
                y = cf0.k
            ] {\orbit};\label{plot:ex3:cf0}

        \addplot [esplot, otherplot, var2] table[
                x = cf1.y,
                y = cf1.k
            ] {\orbit};\label{plot:ex3:cf1}

        \addplot [esplot, otherplot, var3] table[
                x = cf2.y,
                y = cf2.k
            ] {\orbit};\label{plot:ex3:cf2}
    \end{axis}
\end{tikzpicture}

%% file: extended/cf_calculations.tex
\end{tabularx}

%% file: figures/tree_y.tikz
\setlength{\blockheight}{0.28cm}%
\setlength{\blockdistance}{\blockheight + 1pt}%
\setlength{\blockwidth}{1.3\blockheight -1pt}%
%
\tikzstyle{level 1}=[level distance=3cm, sibling distance=27\blockdistance]%
\tikzstyle{level 2}=[level distance=4cm, sibling distance=9\blockdistance]%
\tikzstyle{level 3}=[level distance=4cm, sibling distance=3\blockdistance]%
\tikzstyle{level 4}=[level distance=4cm, sibling distance=\blockdistance]%
%
\tikzstyle{char} = [circle,minimum width=3pt, text centered, font=\tiny, inner sep=1pt]%
\tikzstyle{end} = [circle, minimum width=3pt,fill, inner sep=0pt]%
\tikzstyle{minorderT} = [color=cMRed, font=\tiny]%
\tikzstyle{maxorderT} = [color=cMBlue, font=\tiny]%
\tikzstyle{debug} = [font=\tiny]%
%
%
%
%
\begin{tikzpicture}[%
        grow=right,
        sloped,
        node distance = 0mm and 0mm,
        inner sep = 1pt,
        columnhead/.append style={
                color=cMBlue,
            },
        my shape/.style={
                rectangle split,
                rectangle split horizontal,
                rectangle split parts=#1,
                draw=cLight,
                minimum height=\blockheight,
                rectangle split part align=base,
                every one node part/.style={text width=\blockwidth, align=center, font=\tiny},
                every two node part/.style={text width=\blockwidth, align=center, font=\tiny},
                every three node part/.style={text width=\blockwidth, align=center, font=\tiny},
                every four node part/.style={text width=\blockwidth, align=center, font=\tiny},
                every five node part/.style={text width=\blockwidth, align=center, font=\tiny},
                every six node part/.style={text width=\blockwidth, align=center, font=\tiny},
                fill=white,
            },
        root/.append style={},
        edge from parent/.append style={
                draw = cMedium,
            },
        charting/.append style={
                child anchor=west,
                parent anchor=border
            },
        skipped/.append style={
                color=cLight,
            },
        index/.append style={
                color=black,
            },
        coefficient/.append style={
                color = cNormal,
            },
        powerOfB/.append style={
                color = cNormal,
            },
        powerOfY/.append style={
            },
        powerOfA/.append style={
            },
        powerOfT/.append style={
            }]

        \input{figures/charttree1.tikz}

        \node[columnhead, above=\blockheight] at (node0000.north) (k4) {$d=3$};
        \node[columnhead] at (node000 |- k4) (k3) {$d=2$};
        \node[columnhead] at (node00 |- k4) (k2) {$d=1$};
        \node[columnhead] at (node0 |- k4) (k1) {$d=0$};
        \node[columnhead] at (root |- k4) (k0) {};

        \path (k4) -- coordinate (k43) (k3) -- coordinate (k32) (k2) -- coordinate (k21) (k1) -- coordinate (k10) (k0);
        \draw[dashed, cMedium] (k4 -| k43) -- (k43 |- node2222.south);
        \draw[dashed, cMedium] (k4 -| k32) -- (k32 |- node2222.south);
        \draw[dashed, cMedium] (k4 -| k21) -- (k21 |- node2222.south);
        \draw[dashed, cMedium] (k4 -| k10) -- (k10 |- node2222.south);
%
%
\end{tikzpicture}%

%% file: figures/tree_k.tikz
\setlength{\blockheight}{0.28cm}%
\setlength{\blockdistance}{\blockheight + 1pt}%
\setlength{\blockwidth}{1.5\blockheight -1pt}%
%
\tikzstyle{level 1}=[level distance=3cm, sibling distance=27\blockdistance]%
\tikzstyle{level 2}=[level distance=4cm, sibling distance=9\blockdistance]%
\tikzstyle{level 3}=[level distance=4cm, sibling distance=3\blockdistance]%
\tikzstyle{level 4}=[level distance=4cm, sibling distance=\blockdistance]%
%
\tikzstyle{char} = [circle,minimum width=3pt, text centered, font=\tiny, inner sep=1pt]%
\tikzstyle{end} = [circle, minimum width=3pt,fill, inner sep=0pt]%
\tikzstyle{minorderT} = [color=cMRed, font=\tiny]%
\tikzstyle{maxorderT} = [color=cMBlue, font=\tiny]%
\tikzstyle{debug} = [font=\tiny]%
%
%
%
%
\begin{tikzpicture}[%
        grow=right,
        sloped,
        node distance = 0mm and 0mm,
        inner sep = 1pt,
        columnhead/.append style={
                color=cMBlue,
            },
        my shape/.style={
                rectangle split,
                rectangle split horizontal,
                rectangle split parts=#1,
                draw=cLight,
                minimum height=\blockheight,
                rectangle split part align=base,
                every one node part/.style={text width=\blockwidth, align=center, font=\tiny},
                every two node part/.style={text width=\blockwidth, align=center, font=\tiny},
                every three node part/.style={text width=\blockwidth, align=center, font=\tiny},
                every four node part/.style={text width=\blockwidth, align=center, font=\tiny},
                every five node part/.style={text width=\blockwidth, align=center, font=\tiny},
                every six node part/.style={text width=\blockwidth, align=center, font=\tiny},
                fill=white,
            },
        root/.append style={},
        edge from parent/.append style={
                draw = cMedium,
            },
        charting/.append style={
                child anchor=west,
                parent anchor=border
            },
        skipped/.append style={
                color=cLight,
            },
        index/.append style={
                color=black,
            },
        coefficient/.append style={
                color = cNormal,
            },
        powerOfB/.append style={
                color = cNormal,
            },
        powerOfY/.append style={
            },
        powerOfA/.append style={
            },
        powerOfT/.append style={
            }]

        \input{figures/charttree2.tikz}

        \node[columnhead, above=\blockheight] at (node0000.north) (k4) {$d=3$};
        \node[columnhead] at (node000 |- k4) (k3) {$d=2$};
        \node[columnhead] at (node00 |- k4) (k2) {$d=1$};
        \node[columnhead] at (node0 |- k4) (k1) {$d=0$};
        \node[columnhead] at (root |- k4) (k0) {};

        \path (k4) -- coordinate (k43) (k3) -- coordinate (k32) (k2) -- coordinate (k21) (k1) -- coordinate (k10) (k0);
        \draw[dashed, cMedium] (k4 -| k43) -- (k43 |- node2222.south);
        \draw[dashed, cMedium] (k4 -| k32) -- (k32 |- node2222.south);
        \draw[dashed, cMedium] (k4 -| k21) -- (k21 |- node2222.south);
        \draw[dashed, cMedium] (k4 -| k10) -- (k10 |- node2222.south);
%
%
\end{tikzpicture}